\documentclass [11pt,oneside]{amsart}

\reversemarginpar \textwidth=14cm \textheight=19cm

\usepackage{amssymb} \usepackage{amsfonts} \usepackage{amsmath}
\usepackage{amsthm} \usepackage{epsfig} 
\usepackage[cmtip,arrow]{xy}
\usepackage{amsmath,amssymb,enumerate,pb-diagram,pb-xy}

\usepackage[applemac]{inputenc}
\input xy
\xyoption{all}

\usepackage{caption}

\newcommand{\tto}[1]{\stackrel{#1}{\longrightarrow}}
\newtheorem{lemma}{Lemma}[section]
\newtheorem{teo}[lemma]{Theorem}
\newtheorem{prop}[lemma]{Proposition}
\newtheorem{cor}[lemma]{Corollary}

\theoremstyle{definition}
\newtheorem{defn}[lemma]{Definition}

\newtheorem{rem}[lemma]{Remark}

\newcommand{\calZ}{\ensuremath {\mathcal{Z}}}
\newcommand{\calU}{\ensuremath {\mathcal{U}}}
\newcommand{\matN}{\ensuremath {\mathbb{N}}}
\newcommand{\R} {\ensuremath {\mathbb{R}}}

\newcommand{\matC} {\ensuremath {\mathbb{C}}}

\newcommand{\llim} {\ensuremath {\lim_{n\to\infty}}}
\newcommand{\f} {\ensuremath {\mathcal{F}}}
\newcommand{\sd} {\ensuremath {{\rm sd}}}

\newcommand{\cst} {\ensuremath{C^\ast}}
\newcommand{\climst} {\ensuremath{{C}_b^\ast}}
\newcommand{\ccst} {\ensuremath{C_c^\ast}}
\newcommand{\cclimst} {\ensuremath{{C}_{b,c}^\ast}}

\newcommand{\hc} {\ensuremath{H_c}}

\newcommand{\clim} {\ensuremath{C_b}}
\newcommand{\cc} {\ensuremath{{C}_c}}
\newcommand{\cclim} {\ensuremath{C}_{b,c}}

\newcommand{\hstf} {\ensuremath{\mathcal{H}^\ast}}

\newcommand{\cstf} {\ensuremath{\mathcal{C}^\ast}}

\newcommand{\vare} {\ensuremath{\varepsilon}}

\newcommand{\hst} {\ensuremath{H^\ast}}
\newcommand{\hcst} {\ensuremath{H_c^\ast}}
\newcommand{\hlimst} {\ensuremath{{H}_b^\ast}}
\newcommand{\hclimst}{\ensuremath{{H}_{b,c}^\ast}}

\newcommand{\hcstf} {\ensuremath{\mathcal{H}_c^\ast}}

\newcommand{\ccstf} {\ensuremath{\mathcal{C}_c^\ast}}

\newcommand{\ccf} {\ensuremath{\mathcal{C}_c}}

\newcommand{\hbst} {\ensuremath{H_b^\ast}}

\newcommand{\calC} {\ensuremath {\mathcal{C}}}

\newcommand{\calH} {\ensuremath {\mathcal{H}}}
\newcommand{\hx} {\ensuremath {h_{\widetilde{X}}}}
\newcommand{\xtil} {\ensuremath {\widetilde{X}}}
\newcommand{\res} {\ensuremath {{\rm res}}}

\newcommand{\tr} {\ensuremath {{\rm trans}}}

\author{Roberto Frigerio}

\address{Dipartimento di Matematica \\
Universit\`a di Pisa \\
Largo B.~Pontecorvo 5 \\
56127 Pisa, Italy}

\email{frigerio@dm.unipi.it}

\title[(Bounded) continuous cohomology and proportionality principle]{(Bounded) continuous 
cohomology\\ and Gromov's proportionality principle}

\subjclass[2000]{55N10 (53C23, 55N40, 57N65)}

\keywords{Compact--open topology, continuous cohomology,
Riemannian manifolds, simplicial volume}

\thanks{}

\begin{document}


\begin{abstract}
Let $X$ be a topological space, and let $\cst (X)$ be the complex of singular cochains
on $X$ with coefficients in $\R$. We denote by $\ccst (X)\subseteq \cst (X)$ 
the subcomplex given by \emph{continuous} 
cochains, \emph{i.e.}~by such cochains whose restriction
to the space of simplices (endowed with the compact--open topology) defines a continuous 
real function.
We prove that
at least for ``reasonable'' spaces the inclusion $\ccst (X)\hookrightarrow\cst (X)$
induces an isomorphism in cohomology, thus answering a question
posed by Mostow. 
We also prove that such isomorphism is isometric with respect to the
$L^\infty$--norm on cochains defined by Gromov. 

As an application, we discuss a cohomological
proof of Gromov's proportionality principle for the simplicial volume
of Riemannian manifolds.
\end{abstract}

\maketitle

\section{Preliminaries and statements}\label{construction:section}

Let $X$ be a topological space. We denote by $C_\ast (X)$ (resp.~by $\cst (X)$) the usual
complex of singular chains (resp.~cochains) on $X$ with coefficients in $\R$. 
For $i\in\matN$, we let $S_i (X)$ be the set of singular $i$--simplices
in $X$, and we endow $S_i (X)$ with the compact--open topology (see Appendix~\ref{basi:app} for
basic definitions and results about the compact--open topology). We also
regard $S_i (X)$ as a subset of $C_i (X)$, so that
for any cochain $\varphi\in C^i (X)$ it makes sense to consider its restriction
$\varphi|_{S_i (X)}$. For every $\varphi\in C^i (X)$, we set
$$
|| \varphi ||=||\varphi||_\infty  = \sup \left\{|\varphi (s)|,\, s\in S_i (X)\right\}\in [0,\infty].
$$
We denote by $\climst (X)$ the submodule of bounded cochains, \emph{i.e.}~we set
$\climst (X)= \left\{\varphi\in\cst (X)\, | \, ||\varphi||<\infty\right\}$. Since
the differential takes bounded cochains into bounded cochains, $\climst (X)$
is a subcomplex of $\cst(X)$.
We introduce the following
submodules of $\cst (X)$, which in fact are easily seen to be subcomplexes
of $\cst (X)$:
$$
\begin{array}{cll}
\ccst (X)&=& \left\{\varphi\in\cst (X)\, | \, \varphi|_{S_\ast (X)}\ {\rm is\ continuous}\right\},\\
\cclimst (X)&=& \ccst (X) \cap \climst (X),\\
\end{array}
$$
We denote by $\hst (X)$, $\hcst (X)$, $\hlimst (X)$, $\hclimst (X)$
respectively 
the homology of the complexes $\cst (X)$, $\ccst (X)$, $\climst (X)$, 
$\cclimst (X)$.
Of course, $\hst (X)$ is the usual singular cohomology module of $X$, while $\hlimst (X)$
is the usual bounded cohomology module of $X$. 
Also note that the norm on $C^i (X)$ descends (after 
the suitable restrictions) to a seminorm on each of the modules 
$\hst (X)$, $\hcst (X)$, $\hbst (X)$, $\hclimst (X)$. More precisely, if
$\varphi\in H$ is a class in one of these modules, which is obtained as a quotient
of the corresponding module of cocycles $Z$, then we set 
$$
||\varphi||=\inf \left\{||\psi||,\, \psi\in Z,\, [\psi]=\varphi\ {\rm in}\ H\right\}.
$$
This seminorm may take infinite values on elements in $\hst (X)$, $\hcst (X)$,
and may be a priori null on non-zero elements in $\hcst (X)$,
$\hlimst (X)$, $\hclimst (X)$
(but not on non-zero elements in $\hst (X)$: it is easy to see that a cohomology class with norm equal to zero
has to be null on any cycle, whence null in $\hst (X)\cong (H_\ast (X))^\ast$).

A map $\alpha\colon E\to F$ between seminormed spaces is \emph{norm--decreasing}
if $||\alpha (v)||\leq ||v||$ for every $v\in E$ (this definition makes sense even when
the seminorms considered take $+\infty$ as value).
The natural inclusions
$$
\ccst (X)\hookrightarrow \cst (X), 
\qquad  \cclimst (X)\hookrightarrow\climst (X),
$$ 
induce norm--decreasing maps in cohomology
$$
i^\ast  \colon  \hcst (X)\to \hst (X), \qquad  i_b^\ast  \colon  \hclimst (X)\to\hlimst (X).
$$ 
Moreover, the natural inclusions
$$
\climst (X) \hookrightarrow \cst (X), \qquad \cclimst (X)\hookrightarrow \ccst (X),
$$ 
induce maps in cohomology
$$
c^\ast  \colon  \hlimst (X)\to \hst (X), \qquad c_c^\ast \colon \hclimst (X)\to\hlimst (X).
$$ 
The map $c^\ast$ (resp.~$c_c^\ast$) is called \emph{comparison map}
(resp.~\emph{continuous comparison map}).
By the very definitions, for $\varphi\in\hst (X)$, $\psi\in\hcst (X)$ we have 
$$
\begin{array}{lll}
||\varphi|| & = & \inf \big\{||\varphi_b||,\, \varphi_b\in\hlimst (X),\, c^\ast (\varphi_b)=\varphi\big\},\\ 
||\psi|| & = & \inf \big\{||\psi_b||,\, \psi_b\in\hclimst (X),\, c_c^\ast (\psi_b)=\psi\big\},\\
\end{array}
$$
where such infima are intended to be equal to $\infty$ when taken over the empty set. Moreover, we obviously have $i^\ast\circ c_c^\ast=c^\ast\circ i_{b}^\ast$.

\subsection{Continuous vs.~usual cohomology}\label{cont:subsec}
It is well-known that bounded cohomology and standard cohomology are very different from each other,
while 
Bott stated in~\cite{Bott} that, at least for ``reasonable spaces'', the map $i^\ast$
is an isomorphism. However, Mostow asserted in~\cite[Remark 2 at p.~27]{Mostow} that the natural proof of this fact
seems to raise some difficulties. 

More precisely, it is quite natural to ask whether continuous cohomology
satisfies Eilenberg-Steenrod axioms for cohomology.
In fact, it is not difficult to show that
continuous cohomology verifies 
the so-called ``dimension axiom'' (see Subsection~\ref{dimen:subsec})
and ``homotopy axiom'' (see Subsection~\ref{hom:subsec}).
However, if $Y$ is a subspace of $X$ it is in general not possible
to extend cochains in $\ccst (Y)$ to cochains in $\ccst (X)$, so that it is not clear if a natural long
exact sequence for pairs actually exists in the realm of continuous cohomology.
This difficulty can be overcome either by considering only pairs $(X,Y)$ where
$X$ is metrizable and
$Y$ is closed in $X$, or by exploiting a cone construction, 
as described in~\cite{russo}.
A still harder issue
arises about excision: even if the barycentric subdivision operator consists of
a finite sum (with signs) of continuous self-maps of $S_\ast (X)$,
the number of times a simplex should be subdivided in order to become ``small''
 with respect to a given open cover 
depends in a decisive way on the simplex itself. Thus the dual map of  
(suitable iterations) of the barycentric subdivision 
does not necessarily carry continuous cochains to continuous cochains.
However, in Proposition~\ref{vannobene:prop} we show 
that such dual map takes \emph{locally zero}
continuous cochains into \emph{locally zero} continuous cochains.
Together with some sheaf--theoretic arguments, 
this 
turns out to be sufficient in order to prove the following:

\begin{teo}\label{stand:teo}
Suppose $X$ has the homotopy type of a metrizable and locally contractible
topological space.
Then the map $i^\ast$ is an isomorphism.
\end{teo}

A similar result  
has recently 
been obtained
in~\cite{russo} under stronger assumptions on $X$ (but considering
a larger class of groups of coefficients) in~\cite{russo} .
Mdzinarishvili has shown that continuous cohomology
satisfies in fact the Eilenberg-Steenrod
axioms for cohomology (at least when considering continuous cohomology as a functor
defined on the category of \emph{metric} spaces).
As expected, Mdzinarishvili's argument for showing
that continuous cohomology satisfies the axiom of excision 
is quite subtle (see also Remark~\ref{nuovorem} below).

\smallskip

In order to show that
the isomorphism $i^\ast$ is isometric, 
at least for a large class
of spaces,
in Section~\ref{algebra1:sec}
we describe how the continuous bounded cohomology and the 
ordinary bounded cohomology
of $X$ are related to the bounded
cohomology of the fundamental group of
of $X$. Building on results and techniques developed by Ivanov~\cite{Ivanov}
and Monod~\cite{Monod}, we are then able to prove the following
result:

\begin{teo}\label{lim2:teo}
Suppose $X$ is path connected and has the homotopy type of a countable CW--complex. Then the map
$i_b^\ast$ is surjective and norm--decreasing. Moreover, $i_b^\ast$
admits a right inverse which is an isometric embedding.
\end{teo}

If $X$ satisfies the hypothesis of Theorem~\ref{lim2:teo}, then
it is metrizable, so by Theorem~\ref{stand:teo}
the map $i^\ast\colon \hcst (X)\to \hst(X)$
is an isomorphism. Moreover, for every $\varphi\in \hst (X)$ 
we have

$$
\begin{array}{lll}
||\varphi||&=&\inf \{||\varphi_b||\ |\ \varphi_b\in \hlimst (X),\, c^\ast (\varphi_b)=\varphi\} \\ &=&
\inf \{||\varphi_{b,c}||\ |\ \varphi_{b,c}\in \hclimst (X),\, c^\ast(i_{b}^\ast (\varphi_{b,c}))=\varphi \}\\ &=&
\inf \{||\varphi_{b,c}||\ |\ \varphi_{b,c}\in \hclimst (X),\, i^\ast(c_c^\ast (\varphi_{b,c}))=\varphi \}\\ &=&
\inf \{||\varphi_{b,c}||\ |\ \varphi_{b,c}\in \hclimst (X),\, c_c^\ast (\varphi_{b,c})=(i^\ast)^{-1} (\varphi) \} \\
& = & ||(i^\ast)^{-1} (\varphi)|| ,
\end{array}
$$
where the second equality is due to Theorem~\ref{lim2:teo}. As a consequence
we get the following:

\begin{teo}\label{standiso:teo}
 Suppose $X$ is path connected and has the homotopy type of a countable CW--complex.
Then the map $i^\ast$ is an isometric isomorphism.
\end{teo}

\subsection{The case of aspherical spaces}
When $X$ has contractible universal covering, 
the relation between the (bounded) cohomology of $X$
and the continuous (bounded) cohomology of $X$ is more explicit
than in the general case. 
More precisely, in Section~\ref{main:sec}
we prove
the following:

\begin{teo}\label{lim:teo}
Suppose $X$ is path connected and paracompact with contractible universal covering.
Then the maps $i^\ast$ and
$i_b^\ast$ are isometric isomorphism.
Moreover,
the inverse maps 
$$
\left(i^\ast\right)^{-1}\colon H^\ast (X)\to H_c^\ast (X),\quad \left(i_b^\ast\right)^{-1}\colon H^\ast_b (X)\to H^\ast_{b,c} (X)
$$
can be described by explicit formulae.
\end{teo}

The techniques developed for the proof of Theorem~\ref{lim:teo} can be adapted to 
provide, in dimension one, a more explicit description of the inverse map
$\left(i^{1}\right)^{-1}\colon H^1 (X)\to H^1_c (X)$, even when
$\xtil$ is not contractible. More precisely, in Subsection~\ref{lowdim2:sub}
we prove the following:

\begin{teo}\label{lowdim2:teo}
Suppose that $X$ is paracompact, locally path connected  
and semilocally simply connected. Then the map
$$
i^1\colon H_c^1 (X)\to H^1 (X)
$$ 
is an isomorphism, whose inverse map can be described by an explicit
formula. Moreover, we have ${H}_{b,c}^1 (X)=0\, (={H}_b^1 (X))$.
\end{teo}

In Section~\ref{1dim:sec}
we give examples of path connected (and in one case even simply connected!) ``pathological'' spaces whose first continuous
cohomology module is not isomorphic (through $i^\ast$) to the standard first cohomology module. Such spaces are suitable variations of the
\emph{Hawaiian earring space} and of the \emph{comb space}.

\subsection{Simplicial volume}\label{simpl:sub}
The $L^\infty$--seminorm on cohomology introduced above naturally arises
as dual to the $L^1$--seminorm on homology we are going to describe.
If $X$ is a topological space and $\alpha\in C_i (X)$, we set
$$
||\alpha||=||\alpha||_1 =\sum_{\sigma\in S_i (X)} |a_\sigma|\, ,\quad {\rm where}\ 
\alpha=\sum_{\sigma\in S_i (X)} a_\sigma \sigma.
$$ 
(Note that the sums in the formula above are indeed finite, due to the definition
of singular chain). 
This norm descends to a seminorm on $H_\ast (X)$, which is defined as follows:
if $[\alpha]\in H_i (X)$, then
$$
||[\alpha ]||  =  \inf \{||\beta ||,\, \beta\in C_i (X),\, d\beta=0,\, [\beta ]=[\alpha ] \}.
$$
This seminorm can be null on non-zero elements of $H_\ast (X)$.

If $X$ is a $n$--dimensional connected closed orientable manifold, 
then we denote by $[X]_\mathbb{R}$ 
the image of a generator
of $H_n (X; \mathbb{Z})\cong \mathbb{Z}$ under the change of coefficients map
$H_n (X; \mathbb{Z})\to H_n (X; \R)=H_n (X)$. 
The \emph{simplicial volume} of $X$ is 
$||X||= ||[X]_\R||$. It is easily seen that if $Y$ is the total space of a $d$--sheeted covering
of $X$, then $||Y||=d\cdot ||X||$. Thus, if $X$ is non-orientable, it is reasonable to set
$||X||=||X'|| / 2$, where $X'$ is the double covering of $X$ with orientable total space.

\subsection{The proportionality principle}
Suppose now $X,Y$ are closed connected Riemannian manifolds. 
If $X,Y$ admit a common finite Riemannian covering,
then multiplicativity of Riemannian and simplicial volumes under finite coverings
implies that $||X||/{\rm Vol}(X)=||Y||/{\rm Vol}(Y)$. Gromov's proportionality principle~\cite{Gromov}
ensures that
this equality still holds even when $X,Y$ only share the universal covering:

\begin{teo}[\cite{Gromov,Loh,Bucher}]\label{prop:teo}
Let $X,Y$ be closed Riemannian manifolds with isometric Riemannian universal covering.
Then
$$
\frac{||X||}{{\rm Vol} (X)}=\frac{||Y||}{{\rm Vol}(Y)}.
$$
\end{teo}

A detailed proof of the proportionality principle has been provided
by L\"oh~\cite{Loh} following the ``measure homology'' approach due to Thurston~\cite{Thurston}, 
while the strategy described by Bucher-Karlsson in~\cite{Bucher} 
makes explicit use of bounded cohomology, more in the spirit of
the original argument by Gromov (however, it may be worth mentioning that, in~\cite{Loh}, 
the proof of the fact that measure homology is isometric to
the standard singular homology, which is the key step towards the proportionality 
principle, still relies on results about bounded cohomology).

However, Gromov's approach to the proportionality principle exploits
an averaging process which can be defined
only on sufficiently regular (\emph{e.g.}~bounded Borel measurable)
cochains. Moreover, one needs a regularity result for those 
cochains that are defined by integrating differential forms.
As a consequence, singular (continuous or Borel)
cohomology has to be replaced by \emph{smooth} singular (continuous or Borel)
cohomology, \emph{i.e.}~by the homology of the complex
of (continuous or Borel)
cochains defined on the set of \emph{smooth} simplices, endowed with the 
\emph{$C^1$--topology} (rather than with the compact--open topology).
These technical details seem to raise some difficulties in 
Bucher-Karlsson's proof of the proportionality principle, which slightly relies
on the
expected isomorphism between ordinary and continuous cohomology 
(attributed in~\cite{Bucher} to Bott~\cite{Bott}), 
and does not deal with the fact that the integration of a volume form
does not define a continuous cochain with respect to the compact--open
topology (see Remark~\ref{noncont:rem}).
Building on the results about smooth continuous
cohomology described in Section~\ref{smooth:sec}, in Section~\ref{gromov:sec}
we closely follow the strategy described in~\cite{Bucher} 
for filling in the details in Gromov's original proof 
of the proportionality principle.

\subsection{Plan of the paper} In Section~\ref{first:sec} we deal with the basic properties
of continuous cohomology. In particular, we prove functoriality
and homotopy invariance of continuous cohomology, we compute the
continuous cohomology of the point, and we discuss
barycentric (co)subdivisions
of continuous cochains. The results obtained are then
used in Section~\ref{bor:sec} for showing, via sheaf--theoretic arguments,
that continuous cohomology is canonically isomorphic to standard singular cohomology
for a large class of spaces. In Section~\ref{algebra1:sec} 
we first introduce the needed algebraic notions
for dealing with the seminorms on continuous cohomology and bounded continuous cohomology introduced above: in particular, we describe
the machinery of injective and relative injective strong resolutions developed in~\cite{Ivanov, Monod}.
This machinery is then used in Section~\ref{main:sec} for proving 
Theorems~\ref{lim2:teo},~\ref{lim:teo} and~\ref{lowdim2:teo}. Examples of spaces whose continuous
cohomology is \emph{not} isomorphic to the singular one are given in Section~\ref{1dim:sec}.
As explained above, 
when dealing with smooth manifolds, it is useful to consider the space
of smooth simplices,
endowed with the $C^1$--topology.
In Section~\ref{smooth:sec} we describe how definitions and results for 
continuous cohomology can be adapted to this setting, proving
results that are used in the final section, which is devoted to Gromov's proportionality
principle for the simplicial volume of Riemannian manifolds.

\subsection*{Acknowledgement} 
The author thanks Alessandro Sisto for a clever suggestion
that has lead to a decisive improvement of the statement of Theorem~\ref{stand:teo}.

\section{Basic properties of continuous cohomology}\label{first:sec}
We begin by showing that all the theories we have introduced indeed provide 
homotopy functors from the category of topological spaces to the category of graded real vector spaces.

\subsection{Functoriality}
By Lemma~\ref{basi:top:lemma1},
if $f\colon X\to Y$ is a continuous map, then $f_\ast |_{S_i (X)} \colon S_i (X)\to S_i (Y)$
is continuous, so
its dual map sends continuous 
cochains to continuous cochains,
thus defining a map $f_c^\ast\colon \hcst (Y)\to \hcst (X)$.
Of course, if
$f=\rm {Id}_X$, then $f_c^\ast={\rm Id}_{\hcst (X)}$, 
and if $g\colon Y\to Z$
is continuous, then $(g\circ f)_c^\ast=f_c^\ast \circ g_c^\ast$.
The same results hold true also
for continuous bounded cohomology.

\begin{rem}\label{naturale}
 Let us denote by $i^\ast_X\colon \hcst (X)\to \hst (X)$,
$i^\ast_Y\colon \hcst (Y)\to \hst (Y)$ the maps induced by the inclusion of continuous cochains into the space of singular cochains. 
With the above notations,
it is readily seen that $i^\ast_Y\circ f_c^\ast=f^\ast\circ i^\ast_X$.
This shows that $i^\ast$ provides a \emph{natural transformation}
from the functor $H_c^\ast (\cdot)$ to the functor $H^\ast (\cdot)$.
The analogous result also holds in the bounded case for the tranformation $i_b^\ast$.
\end{rem}

\subsection{The dimension axiom}\label{dimen:subsec}
Suppose $X$ consists of only one point. Then the space $S_n (X)$ consists of only one point (the constant simplex),
so any cochain is automatically continuous, and
continuous cohomology coincides with the usual singular cohomology theory.  We have thus proved the following:
\begin{prop}
Suppose $X$ consists of only one point. Then
$$
\begin{array}{llllll}
H_c^0 (X)&\cong & {H}_{b,c}^0 (X) &\cong &\R, &\\
H^i_c (X)&=& {H}_{b,c}^i (X) &=& 0 &\quad {\rm for}\ i\geq 1.
\end{array}
$$
\qed
\end{prop}
\smallskip

Let now $X=\bigsqcup_{i\in I} X_i$ be the disjoint union of the topological spaces $X_i$, $i\in I$,
and endow $X$ with the disjoint union topology (so $\Omega\subseteq X$ is open in $X$ if and only 
if $\Omega\cap X_i$ is open in $X_i$ for every $i\in I$).
Since each $X_i$ is open in $X$ and the standard simplex is connected, we have 
$C_\ast (X)=\bigoplus_{i\in I} C_\ast (X_i)$, whence $C^\ast (X)=\prod_{i\in I} C^\ast (X_i)$. Moreover,
we have the natural inclusion 
$l_c\colon C_c^\ast (X)\hookrightarrow \prod_{i\in I} C_c^\ast (X_i)$, 
which restrict to the inclusion
${l}_{b,c}\colon \cclimst (X)\hookrightarrow \prod_{i\in I} \cclimst (X_i)$.

Now, since each $X_i$ is open in $X$,
each $S_n (X_i)$ is open in $S_n (X)$. This readily implies that $\varphi\colon S_n (X)\to\R$
is continuous if and only if $\varphi|_{S_n (X_i)}$ is continuous
for every $i\in I$. Therefore $l_c$ is an isomorphism.
Moreover, if $I$ 
is finite, then $\varphi\colon S_n (X)\to\R$
is bounded if and only if $\varphi|_{S_n (X_i)}$ is bounded, so ${l}_{b,c}$
is also an
isomorphism. We can summarize the above discussion in the following:

\begin{prop}\label{connected:prop}
Suppose $X=\bigsqcup_{i\in I} X_i$ is as above. Then 
$\hcst (X)\cong \prod_{i\in I} \hcst (X_i)$. Moreover, 
if $I$ is finite then
$\hclimst (X)\cong \prod_{i\in I} \hclimst (X_i)$.
\qed
\end{prop}

\subsection{Homotopy invariance}\label{hom:subsec}
We now show that all the cohomology theories we have introduced indeed provide homotopy invariants.

\begin{prop}\label{hom:prop}
Let $f,g\colon X\to Y$ be continuous maps with induced morphisms 
$$
f_c^\ast,g_c^\ast \colon \hcst (Y)\to\hcst (X),\qquad 
{f}_{b,c}^\ast,{g}_{b,c}^\ast 
\colon \hclimst (Y)\to\hclimst (X).
$$
If $f,g$ are homotopic, then
$f_c^\ast=g_c^\ast$, and
${f}_{b,c}^\ast={g}_{b,c}^\ast$.
\end{prop}
\noindent {\sc Proof:}
If $H\colon X\times [0,1]\to Y$ is a homotopy between $f$ and $g$,
then there exists an algebraic homotopy $T_\ast\colon C_\ast (X)\to C_{\ast +1} (Y)$ 
between $f_\ast$ and $g_\ast$
which is constructed as follows: for any $s\in S_i (X)$, the chain $T(s)\in C_{i+1} (Y)$ 
is obtained as the image through $H\circ 
(s\times{\rm Id})$ of a suitable fixed subdivision of the prism $\Delta_i\times [0,1]$. Therefore
$T_\ast$ 
takes any simplex $s\in S_i (X)$ to a fixed number of simplices in $S_{i+1} (Y)$, each of which
continuously depends on $s$ (see Lemma~\ref{basi:top:lemma1}). 
Thus the dual homotopy $T^\ast$ takes continuous cochains into continuous 
cochains, and bounded cochains into bounded cochains, 
whence the conclusion.
\qed\smallskip

\subsection{Barycentric subdivisions}\label{bary:subsec}
The aim of this subsection is to show that the barycentric subdivision operator on
singular chains can be suitably dualized in order to provide a well-defined operator
on the space of \emph{locally zero} continuous cochains.
This fact will play a crucial r\^ole in the proof
of Theorem~\ref{stand:teo}.

Suppose $\calU=\{U_i\}_{i\in I}$ is an open cover of $X$. 
We say that a simplex $s\in C_\ast (X)$
is \emph{$\calU$--small} if its image lies in $U_i$ for some $i\in I$, 
we denote by $S_\ast (X)^\calU\subseteq S_\ast (X)$
the space
of $\calU$--small simplices, and by
$C_\ast(X)^\calU$ the subspace of $C_\ast (X)$ generated by $S_\ast (X)^\calU$.

The usual construction of the barycentric subdivision operator (see \emph{e.g.}~\cite[page 56]{Hu})
provides operators
$$
\sd_n\colon C_n (X)\to C_n (X),\quad D_n\colon C_n (X)\to C_{n+1} (X)
$$
such that the following conditions hold:
\begin{itemize}
\item
for every $s\in S_n (X)$, there exists $k\in\matN$ such that $\sd_n^k (s)\in C_n (X)^\calU$;
\item
$\sd_n \circ d_{n+1}=d_n \circ\sd_n$ and 
$d_{n+1}\circ D_n +D_{n-1} \circ d_n =\sd_n - {\rm Id}_{C_n (X)}$ for every $n\geq 0$,
where $d_\ast\colon C_\ast (X)\to C_{\ast -1} (X)$ is the usual differential
(and $D_{-1}$ is intended to be the null operator);
\item
$\sd_\ast (C_\ast (Y))\subseteq C_\ast (Y)$,
$D_\ast (C_\ast (Y))\subseteq C_{\ast +1} (Y)$ for every
subset $Y$ of $X$.
\end{itemize}
Moreover, 
for any given singular simplex $s\colon \Delta_n\to X$,
the value of $\sd_n (s)$ and of $D_n (s)$ only involves sums 
(with signs) of compositions of suitable
restrictions of $s$ with affine parameterizations of convex subsets of $\Delta_n$. Therefore, 
by Lemma~\ref{basi:top:lemma1}
the restrictions of
$\sd_n$ and $D_n$ to  $S_n (X)$ can both be expressed as algebraic
sums of a finite number of continuous
functions. In particular, 
\begin{itemize}
\item
if $\varphi\in C_c^n (X)$, then 
$\varphi\circ \sd_n \in C_c^n (X)$, $\varphi\circ D_{n-1} \in C_c^{n-1} (X)$.
\end{itemize} 

If $s\in C_\ast (X)$ is any simplex, we set 
$\xi^\calU (s)=\min \{n\in \matN\, |\, \sd_\ast^n (s)\in C_\ast (X)^\calU\}$. 
For every $s\in C_n (X)$, let us denote by $s^0,\ldots,s^n$ the \emph{faces} of $s$ (\emph{i.e.}~the
maps obtained by composing with $s$ the affine inclusions 
of $\Delta_{n-1}$ onto the faces of $\Delta_n$).

Following~\cite{Hu}, 
for every $n\in\matN$ we define the homomorphisms
$$
\tau^\calU_n\colon C_n (X)\to C_n (X)^\calU,\qquad
\Omega^\calU_n\colon C_n (X)\to C_{n+1} (X)
$$ 
as the unique linear maps such that for every $s\in S_n (X)$ 
$$
\begin{array}{lll}
\tau^\calU_n (s)&=&\sd_n^{\xi^\calU (s)} (s)
+\sum_{i=0}^n (-1)^i \cdot \left( \sum_{j=\xi^\calU (s^i)}^{\xi^\calU (s)-1} D_{n-1} (\sd_{n-1}^j (s^i))\right), \\
\Omega^\calU_n (s)&=&\sum_{j=0}^{\xi^\calU (s)-1} D_n (\sd_n^j (s)).
\end{array}
$$
It is easily seen
(and shown in~\cite{Hu}) that 
$\tau^\calU_\ast$ is a chain map, and that, if 
$j^\calU_\ast\colon C_\ast (X)^\calU\to C_\ast (X)$ is the natural inclusion,
then
for every $n\geq 0$ we have
$$
j^\calU_n\circ \tau^\calU_n-{\rm Id}_{C_n (X)}=d_{n+1}\circ \Omega^\calU_n+\Omega^\calU_{n-1}\circ d_n,
$$
\emph{i.e.}~$j^\calU_\ast\circ\tau^\calU_\ast$ is chain homotopic to the identity
of $C_\ast (X)$. Since $\tau^\calU_\ast\circ j^\calU_\ast=
{\rm Id}_{C_\ast (X)^\calU}$, this implies
that $\tau^\calU_\ast$, $j^\calU_\ast$ are chain homotopy equivalences. 

Let now $C^n (X)^\calU$ be the dual space of $C_n (X)^\calU$, and endow $C^\ast (X)^\calU$ with the
usual differential. 
We denote by $\tau_\calU^n\colon C^n (X)^\calU \to C^n (X)$, 
$\Omega_\calU^n\colon C^n (X)\to C^{n-1} (X)$,
$j_\calU^n\colon C^n (X)\to C^n (X)^\calU$ the dual maps of $\tau^\calU_n$, $\Omega^\calU_{n-1}$, $j^\calU_n$ 
respectively.

\begin{rem}\label{nuovorem}
We now would like to prove that
the maps just introduced take continuous cochains into continuous
cochains. However, this is in general not true, due to the fact that
the map $\xi^\calU\colon S_n (X)\to\matN$, taking values in a discrete set,
has no hope to be continuous. Therefore, even if
$\sd^\ast$ and $D^\ast$ preserve continuity of cochains,
the operators $\tau_\calU^\ast$ and $\Omega_\calU^\ast$
may not enjoy this property. Therefore, in order to prove an
excision theorem for continuous cohomology, it seems that
one cannot carry out 
the na\"ive strategy of trying to dualize
the barycentric subdivision operator in order to obtain a homotopy equivalence
between the complex of continuous cochains and the complex given
by the restriction of continuous cochains to $C_\ast (X)^\calU$. 
\end{rem}

Let us now set ${_\calU C_c^\ast} (X)=\ccst (X)\cap \ker j_\calU^\ast$,
so that ${_\calU C_c^\ast} (X)$ is the complex of continuous
cochains which vanish on $\calU$--small chains.
Even if the above remark shows that the operators just introduced
cannot preserve continuity of cochains,
for our purposes it will be sufficient to prove that
$\Omega_\calU^\ast$ takes 
${_\calU C_c^\ast} (X)$ into 
${_\calU C_c^\ast} (X)$:


\begin{prop}\label{vannobene:prop}
For every $n\in\matN$ we have 
$$
\Omega_\calU^n ({_\calU C^n_c} (X))\subseteq  {_\calU C^{n-1}_c} (X).
$$
\end{prop}
\noindent {\sc Proof:}
Take $\varphi\in {_\calU C_c^n} (X)$.
We begin by showing that $\Omega_\calU^n (\varphi)$ is continuous at every simplex
$\overline{s}\in S_{n-1} (X)$.
If $\overline{k}=\xi^\calU (\overline{s})$, then
by the very definitions the chain 
$\sd^{\overline k} (\overline s)$ is $\calU$--small.
Moreover, as mentioned above, for every $s\in S_{n-1} (X)$ the chain $\sd_{n-1}^{\overline{k}} (s)$
is the sum (with signs) of a fixed number of 
singular simplices, each of which continuously depends on $s$.
This easily implies that a neighbourhood $W$ of $\overline{s}$ in $S_n (X)$
exists such that for every $s\in W$ the chain
$\sd_{n-1}^{\overline k} (s)$ is $\calU$--small, or equivalently
 $\xi^\calU (s)\leq \overline{k}$.
More in general,
for every $j\in\matN$ the restrictions of
$\sd^j_{n-1}$ and $D_n$ to  $S_{n-1} (X)$ can both be expressed as 
algebraic sums of a finite number of continuous
functions, and this implies that since $\varphi\in C_c^n (X)$, then 
$\varphi\circ D_{n-1}\circ \sd^j_{n-1} \in C_c^n (X)$.

Therefore, in order to prove that $\varphi$ is continuous on $W$, whence at $\overline{s}$,
it is sufficient to show that for every $s\in W$ we have
\begin{equation}\label{definitiva:eq}
\Omega_\calU^n (\varphi) (s)= 
\left(\sum_{j=0}^{\overline{k}-1} \varphi\circ D_{n-1} \circ \sd_{n-1}^j\right) (s), 
\end{equation}
or, equivalently, that
\begin{equation}\label{definitiva2:eq}
0=  \left(\sum_{j=0}^{\overline{k}-1} \varphi ( D_{n-1} (\sd_{n-1}^j (s)))\right)
- \Omega^n_\calU (\varphi) (s)= \varphi \left( \sum_{j=\xi^{\calU}(s)}^{\overline{k}-1} D_{n-1} (\sd_{n-1}^j (s))\right)
\end{equation}
(the point here is that the right--hand side of equation~\eqref{definitiva:eq},
being
the evaluation on $s$ of a \emph{fixed} sum of continuous functions,
defines
a continuous function of $s$).

However, for every $s\in W$ and $j\geq \xi^\calU (s)$ the chain $\sd_{n-1}^j (s)$ is $\calU$--small, and since
$D_{n-1} (C_{n-1} (Y))\subseteq C_n (Y)$ for every subset $Y$ of $X$,
this implies that the chain $\sum_{j=\xi^{\calU}(s)}^{\overline{k}-1} D_{n-1} (\sd_{n-1}^j (s))$
is $\calU$--small. Since $\varphi$ is null on $\calU$--small chains, this
readily gives equality~\eqref{definitiva2:eq}. 

We have thus proved that $\Omega_\calU^n (\varphi)$ is continuous. In order to conclude
we are left to show that $\Omega_\calU^n (\varphi)$ is locally zero. However, if $c\in C_{n-1} (X)$
is $\calU$--small, then $\Omega_{n-1}^\calU (c)=c$, so 
$\Omega^n_\calU (\varphi) (c)=\varphi (\Omega_{n-1}^\calU (c))=\varphi (c)=0$, whence the conclusion.
\qed\smallskip

We now denote by ${_0 \ccst (X)}$ the complex of continuous \emph{locally zero}
cochains, where we say that a continuous cochain is \emph{locally zero}
if it belongs to ${_\calU \ccst (X)}$ for some open covering $\calU$ of $X$.
As a consequence of Proposition~\ref{vannobene:prop} we get the following result,
which will play a crucial r\^ole in the proof of Theorem~\ref{stand:teo}.

\begin{prop}\label{loc0:prop}
We have $H^\ast ({_0 C^\ast_c} (X))=0$.
\end{prop}
\noindent {\sc Proof:}
Take $[\varphi]\in H^n ({_0 C^\ast_c} (X))$, and let $\calU$ be an open cover of $X$
such that $\varphi\in {_\calU C_c^n (X)}$. 
Since $j^n_\calU (\varphi)=0$, we have
\begin{equation}\label{eq}
-\varphi = \tau^n_\calU (j^n_\calU (\varphi))-\varphi = \Omega_\calU^{n+1}(\delta (\varphi))+
\delta(\Omega^n_\calU (\varphi))= \delta (\Omega^n_\calU (\varphi)),
\end{equation}
where we also used that $\delta(\varphi)=0$. 
By Proposition~\ref{vannobene:prop}, therefore, $\varphi$ is the coboundary
of a continuous cochain which belongs to ${_\calU C^{n-1}_c (X)}$, 
whence to ${_0 C^{n-1}_c (X)}$, and this readily implies the conclusion.
\qed\smallskip

\section{Unbounded continuous cohomology of topological spaces}\label{bor:sec}
This section is entirely devoted to the proof of Theorem~\ref{stand:teo}.
After defining sheafified versions of ordinary cohomology and continuous
cohomology, we will prove in Theorem~\ref{iso1:teo} 
that these cohomology theories are
isomorphic to each other
(this fact was already mentioned in~\cite{Mostow}). Since
ordinary cohomology is isomorphic to its sheafified version
(see Theorem~\ref{bredon:teo}), in order to conclude we will 
exploit Proposition~\ref{loc0:prop} for proving
Theorem~\ref{iso2:teo}, which asserts
that also continuous cohomology is isomorphic to its sheafified version.
This last result provides the missing step in Mostow's approach to 
the proof of Theorem~\ref{stand:teo}.

\subsection{Sheaves and presheaves: preliminaries and notations}
We now introduce some notations and recall some basic results about sheaves and presheaves. 
For further reference see \emph{e.g.}~\cite{Bredon}.
Let $X$ be a topological space. If $F$ is a presheaf of real vector spaces
on $X$, and $U\subseteq X$ is open,
we denote by $F(U)$ the sections of $F$ over $U$. As usual, if 
$V\subseteq U$ are open subsets of $X$ and $\varphi\in F(U)$, we denote by
$\varphi |_V$ the restriction of $\varphi$ to $V$. 

A presheaf $F$ is a \emph{sheaf} if 
the following conditions hold:
\begin{itemize}
\item
if $U=\bigcup_{i\in I} U_i$, where each $U_i$ is open, and $\varphi\in F(U)$ is such that
$\varphi |_{U_i}=0$ for every $U_i$, then $\varphi=0$. 
\item
if $U=\bigcup_{i\in I} U_i$, where each $U_i$ is open, and $\varphi_i\in F(U_i)$ are such that
$\varphi_i |_{U_i\cap U_j}=\varphi_j |_{U_i\cap U_j}$ for every $i,j\in I$,
then there exists $\varphi\in F(U)$ such that $\varphi |_{U_i}=\varphi_i$ for every $i\in I$.
\end{itemize}

Let now $F$ be a presheaf of vector spaces on $X$.
We will now describe the classical \emph{sheafification} procedure, which
canonically associates to $F$ a sheaf $\f$.
If $x\in X$, we denote by $F_x$ the stalk of $F$ at $x$; if $x\in U$ and $\varphi\in F(U)$, we denote
by $[\varphi]_x$ the class of $\varphi$ in $F_x$. If $U\subseteq X$ is open, we say that a map
$\psi\colon U\to\bigsqcup_{x\in U} F_x$ is a \emph{continuous section} over $U$ if 
for every $x\in U$ there exist an open set $V_x$ with $x\in V_x\subseteq U$ and an element
$\varphi_x\in F(V_x)$
such that $\psi(y)=[\varphi_x]_y$ for every $y\in V_x$. 
We denote by $\f(U)$ the set of continuous sections over $U$. It is well-known that the association
$U\mapsto \f (U)$ indeed defines a sheaf (restrictions are easily induced by those of $F$). Moreover,
there exists an obvious morphism of presheaves $\rho_F\colon F\to \f$, which is natural
in the following sense: if
$f\colon F\to G$ is a morphism of preshaves, then a unique morphism of sheaves
$\hat{f}\colon \f\to\mathcal{G}$ exists such that $\rho_G\circ f=\hat{f}\circ \rho_F$. Moreover,
if $F$ is already a sheaf, then $\rho_F$ is an isomorphism.

\subsection{The continuous singular cohomology sheaf}\label{outline:subsec}
Remark~\ref{naturale} and Proposition~\ref{hom:prop} ensure that,
in order to prove Theorem~\ref{stand:teo}, we may assume 
without loss of generality that $X$ 
is 
a metrizable locally contractible topological space.
Therefore, this assumption will be taken
for granted from now until the end of this section.

The presheaf of continuous singular
$n$--cochains on $X$ associates to each open set $U\subseteq X$ the real vector space
$C_c^n (U)$, and to every inclusion of open sets 
the obvious restriction of cochains. We will denote such presheaf by
$C_c^n [X]$, and  
$\ccf^n [X]$ will be the sheaf associated to $C_c^n [X]$. 
If $U\subseteq X$ is open, we will denote simply
by $\ccf^n (U)$ (and not by $\ccf^n [X] (U)$) the space of sections of $\ccf^n [X]$
over $U$. The morphism of presheaves 
$\delta\colon C_c^n [X]\to C_c^{n+1} [X]$ induces a morphism between $\ccf^n [X]$ and $\ccf^{n+1}[X]$,
which we will still denote by $\delta$.
We will denote by $\hcstf(X)$ the homology of the complex
$$
0\to\ccf^0 (X)\tto{\delta^0} \ccf^1 (X)\tto{\delta^1}\ldots
\tto{\delta^n} \ccf^{n+1} (X)\tto{\delta^{n+1}}\ldots
$$
\emph{i.e.}~we will set
$\hcstf(X)=H^\ast(\ccstf (X))$. We denote simply by $\rho_c^\ast$ the map 
$\rho_{\ccst[X]}\colon \ccst [X] \to \ccstf [X]$ described
above, and by $\overline{\rho}_c^\ast\colon \hcst(X)\to\hcstf (X)$ the map induced 
in homology by 
the restriction of  $\rho_c^\ast$ to global sections.

The same procedure just described can be applied to standard (\emph{i.e.}~not necessarily continuous)
cochains: to the preasheaf $C^n [X]$ there is associated the sheaf $\mathcal{C}^n [X]$, 
and the usual differential
on cochains induces a differential on the complex $\cstf(X)$, whose homology will be denoted
by $\hstf(X)$. As in the case of continuous cochains, we have a natural map 
$\overline{\rho}^\ast\colon \hst(X)\to\hstf (X)$. Moreover, 
the mentioned naturality of the sheafification process
provides a chain map of complexes of sheaves 
$\ccstf(X)\to\cstf(X)$, and the restriction of this map to 
global sections
induces
in homology a map $i_{\rm sh}^\ast \colon \hcstf(X)\to\hstf(X)$ 
which makes the following diagram commute:
$$
\xymatrix{
\hcst(X) \ar[r]^{\overline{\rho}_c^\ast} \ar[d]^{i^\ast} & \hcstf(X) \ar[d]^{i_{\rm sh}^\ast}\\
\hst(X) \ar[r]^{\overline{\rho}^{\ast}} & \hstf(X).
}
$$
Observe now that, since $X$ is metrizable, 
the family of all closed subsets of $X$ is \emph{paracompactifying}
in the sense of~\cite[page 21]{Bredon}. As a consequence 
we get the 
following classical result (see \emph{e.g.}~\cite[page 26]{Bredon}): 

\begin{teo}\label{bredon:teo}
The map $\overline{\rho}^\ast\colon H^\ast (X)\to \hstf (X)$ is an isomorphism.
\end{teo}
\qed\smallskip

In order to conclude the proof of Theorem~\ref{stand:teo}, it is now sufficient to prove the following
results:

\begin{teo}\label{iso1:teo}
The map $i_{\rm sh}^\ast \colon \hcstf (X)\to \hstf (X)$ is an isomorphism.
\end{teo}

\begin{teo}\label{iso2:teo}
The map $\overline{\rho}_c^\ast\colon H_c^\ast (X)\to \hcstf (X)$ is an isomorphism.
\end{teo}

\subsection{Proof of Theorem~\ref{iso1:teo}}
We refer to~\cite[Section II.9]{Bredon} for the definition of \emph{soft} and \emph{fine}
sheaf.
Let $R[X]$ be the sheaf of real continuous functions over $X$.
For every open subset $U\subseteq X$, the space $C_c^n (U)$ is naturally
a module over the ring $R(U)$, with the operation
defined by $(f\cdot \varphi) (s) = f(s(e_0))\cdot \varphi (s)$
for every $f\in R(U)$, $\varphi\in C_c^n (U)$, $s\in S_n (U)$, where
$e_0$ is the first vertex of the standard $n$--simplex $\Delta_n$.
So
$C_c^n [X]$ is a presheaf of modules over the sheaf of rings with unit $R[X]$,
and this readily implies that $\mathcal{C}_c^n [X]$ is a sheaf of modules over $R[X]$. 
Now, since $X$ is paracompact, the sheaf $R[X]$ is soft (see~\cite[Example~II.9.4]{Bredon}),
so  $\ccstf [X]$ is fine~\cite[Theorem~II.9.16]{Bredon}, whence acyclic.
The very same argument also applies to $\calC^n [X]$, which therefore is
also acyclic. 

Let now $\widetilde{\R}$ be the constant sheaf on $X$ with stalks
isomorphic to $\R$. We recall from Subsection~\ref{outline:subsec} that the inclusion induces   
a chain map of complexes of sheaves: 
$$
\xymatrix{
0 \ar[r] & \widetilde{\R} \ar[d]^{\rm Id} \ar[r]^{\delta^{-1}} & 
\mathcal{C}^0_b [X] \ar[d]  \ar[r]^{\delta^0} & 
\mathcal{C}^1_b [X] \ar[d] \ar[r]^{\delta^1} & 
\mathcal{C}^2_b [X] \ar[d] \ar[r]^{\delta^2} & \ldots
\\
0 \ar[r] & \widetilde{\R} \ar[r]^{\delta^{-1}} & \mathcal{C}^0 [X] \ar[r]^{\delta^0} & 
\mathcal{C}^1 [X] \ar[r]^{\delta^1} & \mathcal{C}^2 [X] \ar[r]^{\delta^2} &\ldots
}
$$
where $\delta^{-1}$ is induced by the usual augmentation map on presheaves which sends
$t\in \R$ to the $0$--cochain which takes the value $t$ on every $0$--simplex.

Recall now that both standard singular cohomology and continuous
cohomology are homotopy invariant: since 
$X$ is locally contractible, this implies that the rows in the diagram above provide
exact sequences of sheaves. Moreover, we have seen that both $\calC^n_c [X]$ and
$\calC^n [X]$ are acyclic for every $n\geq 0$. 
A classical result of sheaf theory
(see \emph{e.g.}~\cite[Theorem II.4.1]{Bredon}) now ensures
that $i_{\rm sh}^\ast\colon \hcstf (X)\to \hstf (X)$ is an isomorphism.
\qed\smallskip

\subsection{Proof of Theorem~\ref{iso2:teo}}
We begin with the following:
\begin{lemma}\label{borelext:lemma}
Let $\calZ$ be a locally finite closed cover of $X$, and 
suppose we are given an element
$\varphi_Z\in C^n_c (Z)$ for every $Z\in \calZ$, in such a way that 
$$\varphi_Z |_{Z\cap Z'}=\varphi_{Z'} |_{Z\cap Z'}\quad {\rm for\ every}\ Z,Z'\in\calZ.$$
Then there exists $\varphi\in C^\ast_c (X)$
such that $\varphi|_{Z} =\varphi_Z$ for every $Z\in \calZ$. 
\end{lemma}
\noindent {\sc Proof:}
We prove the following equivalent statement:
let $\{f_Z\colon S_n (Z)\to\R,\, Z\in\calZ\}$ be a family of continuous functions such that $f_Z|_{S_n (Z\cap Z')}= f_{Z'} |_{S_n (Z\cap Z')}$
for every $Z,Z'\in\calZ$; then 
a continuous function $f\colon S_n (X)\to\R$ exists such that
$f|_{S_n (Z)}=f_Z$ for every $Z\in\calZ$.

Now, since $\calZ$ is a locally finite closed cover of $X$, the family
$\{S_n (Z)\, |\, Z\in\calZ\}$ provides a locally finite collection of
closed subsets of $S_n (X)$. As a consequence, the set 
$W=\bigcup_{Z\in\calZ} S_n (Z)$ is closed in $S_n (X)$, and there exists
a well-defined continuous map
$g\colon W\to \R$ such that $g|_{S_n (Z)}=f_Z$ for
every $Z\in \calZ$. Since $X$ is metrizable and the standard $n$--simplex
is compact, the space $S_n (X)$
is metrizable, so Tietze's extension theorem~\cite[p.~149]{Dug}
ensures that $g$ can be extended to a continuous function defined on the whole
of $S_n (X)$.
\qed\smallskip

\begin{lemma}\label{surj:lemma}
The map $\rho^\ast_c\colon C_c^\ast (X)\to\ccstf (X)$ is surjective.
\end{lemma}
\noindent {\sc Proof:}
Take $\varphi\in \ccf^n (X)$. By the very definitions, for every $x\in X$ an open set $U_x\ni x$ in $X$
and a section $\psi_x\in C_c^n (U_x)$ exist such that $\varphi|_{U_x}=\rho^n_c (\psi_x)$. 
Since $X$ is metrizable, it is paracompact, so there exist
a locally finite open covering $\{V_i\}_{i\in I}$ of $X$ and a function $x\colon I\to X$
such that $\overline{V}_i\subseteq U_{x(i)}$. For every $y\in X$ we set $I(y)=\{i\in I\, |\, 
y\in \overline{V}_i\}$. Since $\{V_i\}_{i\in I}$ is locally finite also the family 
$\{\overline{V}_i\}_{i\in I}$ is locally finite, so 
the set $I(y)$ is finite. Moreover the set $\bigcup_{i\notin I(y)} \overline{V}_i$
is closed, so 
an open neighbourhood $W_y$ 
of $y$ exists such that $W_y\cap\overline{V}_i=\emptyset$ for every $i\notin I(y)$, and this readily implies
that
$I(y')\subseteq I(y)$ for every $y'\in W_y$.
Since $\overline{V}_i\subseteq U_{x(i)}$ for every $i\in I$, up to shrinking $W_y$ we may also assume that
$W_y\subseteq U_{x(i)}$ for every $i\in I(y)$ and
$\psi_{x(i)}|_{W_y}=\psi_{x(j)}|_{W_y}$ for every $i,j\in I(y)$.
For every $y\in X$ we now set $\psi'_y=\psi_{x(i)}|_{W_y}\in C_c^n (W_y)$ for some
$i\in I(y)$ (by construction, this definition does not depend on $i\in I(y)$). 
We now claim that the collection of sections $\psi'_y\in C^n_c (W_y)$, $y\in X$ satisfies the following 
properties:
\begin{enumerate}
\item
$\rho^n_c (\psi'_y)=\varphi|_{W_y}$ for every $y\in X$;
\item
$\psi'_y|_{W_y\cap W_{y'}}=\psi'_{y'}|_{W_y\cap W_{y'}}$ for every $y,y'\in X$
such that $W_y\cap W_{y'}\neq\emptyset$.
\end{enumerate}
Property~(1) is obvious, while if $\emptyset\neq W_y\cap W_{y'}\ni z$, then 
$I(y)\cap I(y')\supseteq I(z)\neq\emptyset$: in particular, if $i_0\in I(y)\cap I(y')$, then
$\psi'_y|_{W_y\cap W_{y'}}= \psi_{x(i_0)}|_{W_y\cap W_{y'}}=\psi'_{y'}|_{W_y\cap W_{y'}}$,
whence property~(2). 

Let now $\mathcal{Z}=\{Z_j\}_{j\in J}$ be a locally finite open cover of $X$
such that for every $j\in J$ there exists $y(j)\in X$ such that $\overline{Z}_j\subseteq W_{y(j)}$.
By property~(2) above and
Lemma~\ref{borelext:lemma}, 
a global cochain $\psi''\in C_c^n (X)$ 
exists such that $\psi''|_{Z_j}=\psi'_{y(j)}|_{Z_j}$ for every $j\in J$.
Now, by property~(1) above, for every $j\in J$ 
we have $\rho^n_c (\psi'')|_{Z_j}=\rho^n_c (\psi''|_{Z_j})=\varphi|_{Z_j}$.
Since $\mathcal{C}^n_c [X]$ is a sheaf, this readily implies that
$\rho^n_c (\psi'')=\varphi$, whence the conclusion.
\qed\smallskip 

We can now conclude the proof of Theorem~\ref{iso2:teo}.
We defined in Subsection~\ref{bary:subsec} the subcomplex ${_0 C^\ast_c (X)}$ of 
\emph{locally zero} continuous cochains.
By its very definition, ${_0 C^\ast_c} (X)$ is equal to the kernel of the map
$\rho^\ast_c\colon C^\ast_c (X)\to \ccstf (X)$. By Lemma~\ref{surj:lemma}
the short sequence of complexes
$$
\xymatrix{
0 \ar[r]& {_0 C^\ast_c} (X) \ar@{^{(}->}[r] & C^\ast_c (X) \ar[r]^{\rho^\ast_c} & \ccstf (X) \ar[r] & 0
}
$$
is therefore exact. This gives a long exact sequence
$$
\xymatrix{
\ldots \ar[r] 
& H^n ({_0 C^\ast_c} (X)) \ar[r] & 
H^n_b (X) \ar[r]^{\overline{\rho}^n_c} & \calH^n_c (X) \ar[r] & H^{n+1} ({_0 C^\ast_c} (X))\ar[r] 
& \ldots
}
$$
By Proposition~\ref{loc0:prop}, we now have $H^n ({_0 C^\ast_c} (X))=H^{n+1} ({_0 C^\ast_c} (X))=0$,
so $\overline{\rho}^n_c$ is an isomorphism.
\qed\smallskip

\begin{rem}
The hypothesis that $X$ is metrizable 
came into play only in the proof of Lemma~\ref{borelext:lemma},
where, in order to extend continuous cocycles defined
on closed subspaces of $S_n (X)$, we exploited the fact that
the space $S_n (X)$ is normal. One could wonder if 
normality of $S_n (X)$  could be proved under
somewhat weaker hypotheses on $X$. However, it is exhibited in~\cite{controesempio}
an example of a paracompact locally contractible space $X$ such that
$S_1 (X)$ is not normal. This seems to suggest that metrizability is 
the most reasonable assumption on $X$ which ensures that $S_n (X)$
is normal for every $n\in\matN$.
\end{rem}

\section{Bounded continuous cohomology of topological spaces}\label{algebra1:sec}
This section is mainly devoted to the proofs of Theorems~\ref{lim2:teo},~\ref{lim:teo}.
We begin by 
reviewing some definitions introduced in~\cite{Ivanov, Monod}.

Let $G$ be a group (which should be thought as endowed with the discrete topology).
In what follows, a \emph{$G$--module} (resp.~a \emph{Banach $G$--module}) is a real vector space (resp.~a real
Banach space) endowed with an action of $G$ (by isometries, in the Banach case) on the \emph{left}.
Sometimes, we will stress the fact that a $G$--module $E$ is \emph{not} endowed with a Banach structure
by saying that $E$ is an \emph{unbounded} $G$--module.

If $E$ is a (Banach) $G$--module, we denote by $E^G\subseteq E$ the submodule
of $G$--invariant elements in $E$, \emph{i.e.}~we set
$$
E^G=\{v\in E\, | \, g\cdot v=v\ {\rm for\ every}\ g\in G\}.
$$  
A \emph{$G$--map} between (Banach) $G$--modules is a (bounded) $G$--equivariant linear map.

\subsection{Relative injectivity}
A $G$--map $\iota\colon A\to B$ between Banach $G$--modules is \emph{strongly injective}
if there exists a linear map $\sigma\colon B\to A$ with $||\sigma||\leq 1$
such that $\sigma\circ\iota={\rm Id}_A$ (we do not require $\sigma$ to be a $G$--map; note
also that strongly injective obviously implies injective).
We now define the important notion of relative injectivity (resp.~injectivity) 
for Banach (resp.~for unbounded)
$G$--modules. 

\begin{defn}\label{inj:def}
A Banach $G$--module $U$ is \emph{relatively injective} if the following holds: 
whenever $A,B$ are Banach $G$--modules, $\iota\colon A\to B$ is a strongly injective $G$--map
and $\alpha\colon A\to U$ is a $G$--map, there exists a $G$--map $\beta\colon B\to U$
such that $\beta\circ\iota = \alpha$ and $||\beta ||\leq ||\alpha||$.
\end{defn}

$$
\xymatrix{
0 \ar[r] & A  \ar[r]_{\iota} \ar[d]_{\alpha} & B \ar@{-->}[dl]^{\beta} \ar@/_/[l]_{\sigma}\\
& U
}
$$

\begin{defn}\label{inj2:def}
An unbounded $G$--module $U$ is \emph{injective} if the following holds: 
whenever $A,B$ are unbounded $G$--modules, $\iota\colon A\to B$ is an injective $G$--map
and $\alpha\colon A\to U$ is a $G$--map, there exists a $G$--map $\beta\colon B\to U$
such that $\beta\circ\iota = \alpha$.
\end{defn}

Note that any injective map between unbounded $G$--modules admits a (maybe not $G$--equivariant)
left inverse, so the notion of relative injectivity can be considered an extension
to the Banach setting of the notion of injectivity for unbounded modules.

\subsection{Resolutions}
A \emph{(Banach) $G$--complex} (or simply a \emph{(Banach) complex}) is a sequence of (Banach)
$G$--modules $E^i$ and $G$--maps $\delta^i\colon
E^i\to E^{i+1}$ such that $\delta^{i+1}\circ\delta^i=0$ for every $i$, where $i$ runs over
$\matN\cup\{-1\}$:

$$
0\longrightarrow E^{-1}\tto{\delta^{-1}} E^0\tto{\delta^0} E^1\tto{\delta^1}\ldots
\tto{\delta^n} E^{n+1}\tto{\delta^{n+1}}\ldots
$$
Such a sequence will be often denoted by $(E^\ast,\delta^\ast)$.

A \emph{chain map} between (Banach) $G$--complexes $(E^\ast,\delta_E^\ast)$ and $(F^\ast,\delta_F^\ast)$
is a sequence of $G$--maps $\{\alpha^i\colon E^i\to F^i,\, i\geq -1\}$
such that $\delta_F^i\circ\alpha^i=\alpha^{i+1}\circ\delta_E^{i}$ for every $i\geq -1$. 
If $\alpha^\ast,\beta^\ast$ are chain maps between $(E^\ast,\delta_E^\ast)$ and $(F^\ast,\delta_F^\ast)$ which coincide in degree $-1$, a 
\emph{$G$--homotopy} 
between $\alpha^\ast$ and $\beta^\ast$ is 
a sequence of $G$--maps $\{T^i \colon E^i\to F^{i-1},\, i\geq 0\}$ such that
$\delta_F^{i-1}\circ T^i + T^{i+1}\circ \delta_E^i=\alpha^i-\beta^i$   for every $i\geq 0$,
and $T_0\circ \delta_E^{-1}=0$. We recall that, according to our definition
of $G$--maps for Banach modules, both chain maps between Banach $G$--complexes
and $G$--homotopies between such chain maps have to be bounded (more precisely, such maps have to be bounded in every degree, while there does not need to be a
uniform bound on their norms as maps from $\oplus_{i\geq -1} E^i$ to $\oplus_{i\geq -1} F^i$).

A complex is \emph{exact} if $\delta^{-1}$ is injective and 
$\ker \delta^{i+1}={\rm Im}\, \delta^i$ for every $i\geq -1$.
Let $E$ be a (Banach) $G$--module. A \emph{resolution} of $E$ as a (Banach) $G$--module
is an exact (Banach) $G$--complex $(E^\ast,\delta^\ast)$ with $E^{-1}=E$.

A resolution 
$(E^\ast,\delta^\ast)$ is \emph{relatively injective} (resp.~\emph{injective})
if $E^n$ is relatively injective (resp.~\emph{injective}) for every $n\geq 0$.

A \emph{contracting homotopy} for a resolution $(E^\ast,\delta^\ast)$ is a 
sequence of linear maps $k^i\colon E^i\to E^{i-1}$ such that $\delta^{i-1}\circ k^i+
k^{i+1}\circ\delta^i = {\rm Id}_{E^i}$ if $i\geq 0$, and $k_0 \circ \delta^{-1}=
{\rm Id}_E$. If $(E^\ast,\delta^\ast)$ is a resolution of Banach modules, the condition
$||k^i||\leq 1$ is also required. 

$$
\xymatrix{
0\ar[r] & E^{-1}\ar[r]_{\delta^{-1}} & E^0 \ar[r]_{\delta^0} \ar@/_/[l]_{k_0} &
E^1  \ar[r]_{\delta^1} \ar@/_/[l]_{k_1} & \ldots \ar@/_/[l]_{k_2} \ar[r]_{\delta^{n-1}}&
E^{n}
\ar[r]_{\delta^n} \ar@/_/[l]_{k_n} & \ldots \ar@/_/[l]_{k_{n+1}}
}
$$

Note however that it is not required that 
$k_i$ is $G$--equivariant. A resolution of a (Banach) $G$--module is \emph{strong}
if it admits a contracting homotopy.

The following results can be proved by means of standard homological algebra arguments
(see~\cite{Ivanov}, \cite[Lemmas 7.2.4 and 7.2.6]{Monod} for the details in the case of Banach resolutions).

\begin{prop}\label{ext:prop}
Let $\alpha\colon E\to F$ be a $G$--map between Banach (resp.~unbounded)
$G$--modules,
let $(E^\ast,\delta_E^\ast)$ be a strong resolution of $E$, and suppose
$(F^\ast,\delta_F^\ast)$ is a $G$--complex with $F^{-1}=F$ and $F^i$ relatively injective
(resp.~injective) for evey $i\geq 0$. 
Then $\alpha$ extends to a chain map $\alpha^\ast$, and any two extensions
of $\alpha$ to chain maps are $G$--homotopic.
\end{prop}

\subsection{(Bounded) group cohomology}
We recall that if $E$ is a (Banach) $G$--module, we denote by $E^G\subseteq E$ the submodule
of $G$--invariant elements in $E$.

Let $(E^\ast,\delta^\ast)$ be a relatively injective strong resolution of the
trivial Banach $G$--module $\R$ (such a resolution exists, see Subsection~\ref{standard:subsec}).
Since coboundary maps are $G$--maps, they restrict to the $G$--invariant submodules of the
$E^i$'s. Thus $((E^\ast)^G,\delta^\ast |)$ is a subcomplex of $(E^\ast,\delta^\ast)$. 
A standard application of Proposition~\ref{ext:prop} now shows that the isomorphism type of
the homology
of $((E^\ast)^G,\delta^\ast|)$ does not depend on the chosen resolution (while the seminorm induced on 
such homology module by the norms on the $E^i$'s could depend on it, see Proposition~\ref{norm:prop} below).
For every $i\geq 0$, we now define the $i$--dimensional \emph{bounded cohomology}
module $H_b^i (G)$ of $G$ (with real coefficients) as follows: if $i\geq 1$, then 
$H_b^i (G)$ is the $i$--th homology module of the complex $((E^\ast)^G,\delta^\ast |)$, while if $i=0$ then $H_b^i (G)=\ker \delta^0\cong \R$.

The same construction applies \emph{verbatim} when considering an injective strong resolution
$(E^\ast,\delta^\ast)$
of $\R$ as an unbounded $G$--module.  In this case, the homology of $((E^\ast)^G,\delta^\ast |)$
is the standard cohomology of $G$, and will be denoted by $\hst (G)$.

\subsection{The standard $G$--resolutions}\label{standard:subsec}
For every $n\in\matN$, let $F^n (G)=\{f\, |\, f\colon G^{n+1}\to\R\}$ and
${F}_b^n (G)=\{f\in F^n (G)\, | \, f\ {\rm is\ bounded}\}$, and endow ${F}_b^n (G)$ with
the supremum norm, thus obtaining a real Banach space. Let $G$ act on $F^n (G)$ in such a way
that $(g\cdot f)(g_0,\ldots,g_n)=f(g^{-1} g_0,\ldots, g^{-1}g_n)$. It is easily seen
that this action leaves ${F}_b^n (G)$ invariant, and endows $F^n (G)$ (resp.~${F}_b^n (G)$)
with a structure of $G$--module (resp.~of Banach $G$--module).
For $n\geq 0$, define $\delta^n \colon F^n (G)\to F^{n+1} (G)$ by setting:
$$
\delta^n (f) (g_0,g_1,\ldots,g_{n+1})=\sum_{i=0}^{n+1} (-1)^i f(g_0,\ldots,\widehat{g}_i,\ldots,g_{n+1}).
$$
It is easily seen that $\delta ({F}_b^n(G))\subseteq {F}_b^{n+1}(G)$, so it makes sense to define
${\delta}_b^n\colon {F}_b^n (G)\to {F}_b^{n+1} (G)$ by restricting $\delta^n$. 
Moreover, we let
$F^{-1} (G)=\R$
be a trivial unbounded $G$--module and ${F}_b^{-1}(G)=\R$ be a trivial Banach $G$--module, and 
we define ${\delta}_b^{-1}\colon \R\to {F}_b^0 (G)$ by
setting ${\delta}_b^{-1} (t) (g)=t$ for every $g\in G$, and $\delta^{-1}$ by composing
${\delta}_b^{-1}$ with the inclusion of ${F}_b^0 (G)$ in $F^0 (G)$.

\begin{rem}
With slightly different conventions and notations, Ivanov proved
in~\cite{Ivanov} that the complex
$({F}_b^\ast (G),{\delta}_b^\ast)$ provides a relatively injective strong
resolution of $\R$ as a bounded $G$--module. His argument can be easily
adapted for showing that the complex
$(F^\ast (G),\delta^\ast)$
is an injective strong resolution of $\R$
as an unbounded $G$--module.
However, these results will not be necessary for our purposes.
\end{rem}


The resolution $(F^\ast (G),\delta^\ast)$ 
(resp.~$({F}_b^\ast (G),{\delta}_b^\ast)$) is usually known
as the standard $G$--resolution of $\R$ as an unbounded $G$--module
(resp.~as a Banach $G$--module).
The seminorm induced on $\hlimst (G)$ by the standard bounded resolution is known as the
\emph{canonical seminorm}.
The following result~\cite{Ivanov, Monod} 
gives a useful characterization of the canonical
seminorm, and plays a decisive r\^ole in our
proof of Theorem~\ref{lim2:teo}.

\begin{prop}\label{norm:prop}
Let $(E^\ast,\delta^\ast)$ be any strong resolution of $\R$
as a Banach $G$--module. Then the identity of $\R$ can be extended to a chain map 
$\alpha_b^\ast$ between $E^\ast$ and the standard resolution of $\R$ as
Banach $G$--module,
in such a way that $||\alpha_b^n||\leq 1$ for every $n\geq 0$. In particular,
the canonical seminorm is not bigger than the seminorm induced on $\hlimst (G)$ by
any
relatively injective strong resolutions.
\end{prop}
\noindent{\sc Proof:}
One can define $\alpha_b^n$ by induction setting, for every $v\in E^n$ and $g_j\in G$:  
$$
\alpha_b^n (v) (g_0,\ldots,g_n)=\alpha_b^{n-1}(g_0(k^{n}(g_0^{-1}(v))))(g_1,\ldots,g_n),
$$
where $\{k^n\}_{n\in\matN}$ is a contracting homotopy for $E^\ast$. It is not difficult 
to prove by induction that  
$\alpha^\ast$ is indeed a norm--decreasing chain $G$--map (see~\cite{Ivanov}, \cite[Theorem 7.3.1]{Monod}
for the details).
\qed\smallskip

\begin{rem}\label{unbounded:rem}
It is readily seen that,
if $(E^\ast,\delta^\ast)$ is any strong resolution of $\R$
as an unbounded $G$--module,  then the formula described in the proof of Proposition~\ref{norm:prop}
also provides an extension of the identity of $\R$ to a chain map
$\alpha^\ast$ between 
$E^\ast$ and the standard resolution of $\R$ as
unbounded $G$--module. Moreover,
$\alpha^n$ is norm--decreasing for every $n\geq 0$. 
\end{rem}

\subsection{Some notations and a useful lemma}\label{notation:subsec}
From now until the end of the section, we assume that $X$ is a path connected paracompact
topological space with universal covering
$p\colon\widetilde{X}\to X$, and we fix an identification between
the fundamental group of $X$ and
the group $\Gamma$ of covering automorphisms of $\widetilde{X}$. Thus every $g\in\Gamma$ defines a 
chain map $g_\ast\colon C_\ast (\xtil)\to C_\ast (\xtil)$.
It is a standard fact of algebraic topology that the action of $\Gamma$ on $\xtil$ is \emph{wandering},
\emph{i.e.}~any $x\in\xtil$ admits a neighbourhood $U_x$ such that $g(U_x)\cap U_x=\emptyset$
for every $g\in\Gamma\setminus \{1\}$ (if $X$, whence $\widetilde{X}$, is locally
compact, then an action on $\widetilde{X}$ is wandering if and only
if it is free and proper). In the following lemma we describe a particular
instance of \emph{generalized Bruhat function} (see~\cite[Lemma 4.5.4]{Monod} for a more general
result based on~\cite[Proposition 8 in VII $\S$2 N$^\circ$ 4]{bou}).

\begin{lemma}\label{tecnico}
There exists a continuous map 
$h_{\widetilde{X}}\colon \xtil\to [0,1]$ with the following properties:
\begin{enumerate}
\item
For every $x\in\xtil$ there exists a neighbourhood $W_x$ of $x$ in $\xtil$
such that the set
$
\left\{g\in\Gamma\, | \,  g(W_x) \cap
{\rm supp}\, h_{\widetilde{X}}
\neq\emptyset\right\}
$
is finite.
\item
For every $x\in\xtil$, we have 
$
\sum_{g\in \Gamma} h_{\widetilde{X}} (g\cdot x)=1
$
(note that the sum on the left-hand side is finite by (1)).
\end{enumerate}
\end{lemma}
\noindent {\sc Proof:}
Let us take a locally finite open 
cover $\{U_i\}_{i\in I}$ of $X$ such that
for every $i\in I$ 
there exists $V_i\subseteq \xtil$ with $p^{-1} (U_i)=\bigcup_{g\in \Gamma}
g(V_i)$ and $g(V_i)\cap g'(V_i)=\emptyset$ whenever $g\neq g'$.
Let $\{\varphi_i\}_{i\in I}$ be a partition of unity adapted to $\{U_i\}_{i\in I}$. It is easily seen that
the map $\psi_i\colon \widetilde{X}\to\R$ which concides with $\varphi_i\circ p$ on $V_i$
and is null elsewhere is continuous. We can now set $h_{\widetilde{X}}=\sum_{i\in I} \psi_i$.
Since $\{U_i\}_{i\in I}$ is locally finite, also $\{V_i\}_{i\in I}$,
whence $\{ {\rm supp}\, \psi_i\}_{i\in I}$, is locally finite, so $h_{\widetilde{X}}$
is indeed well-defined and continuous. 

In order to show that $h_{\widetilde{X}}$ satisfies $(1)$, let $x\in \xtil$, and suppose
$p(x)\in U_{i_0}$. Then there exists $g_0 \in\Gamma$ such that $x\in g_0(V_{i_0})$.
We set $W_x=g_0(V_{i_0})$, and let $J=\{j\in I\, | \, U_j\cap U_{i_0}\neq \emptyset\}$. By construction,
$J$ is finite. Now, if $i\in I\setminus J$ then for every $g\in\Gamma$ we have 
$
p(g(W_x)\cap V_i) \subseteq U_{i_0}\cap U_i=\emptyset
$,
so if $g(W_x)\cap {\rm supp}\, h_{\widetilde{X}}\neq\emptyset$ then $g (W_x)\cap V_{j}
\neq \emptyset$
for some $j\in J$. However, since $g(W_x)\cap g' (W_x)=\emptyset$ for every $g\neq g'$,
for every $\overline{j}\in J$ there is at most one $g\in\Gamma$ such that $g (W_x)
\cap V_{\overline{j}}\neq \emptyset$,
so $\left\{g\in\Gamma\, | \, g(W_x) \cap {\rm supp}\, h_{\widetilde{X}}\neq\emptyset\right\}$
is finite.

Finally, for every $x\in\xtil$, $i\in I$ we have by construction $\sum_{g\in\Gamma} \psi_i (g(x))=
\varphi_i (p(x))$,
so 
$$
\sum_{g\in\Gamma} h_{\widetilde{X}} (g(x))=\sum_{g\in\Gamma} \left(\sum_{i\in I} \psi_i (g(x))\right)=
\sum_{i\in I} \left( \sum_{g\in\Gamma} \psi_i (g(x))\right)=\sum_{i\in I} \varphi_i (p(x))=1,
$$
whence~(2). 
\qed\smallskip

\subsection{Singular cochains as (relatively) injective modules}
For every $n\geq 0$, we define an action of $\Gamma$ 
on $C^n (\xtil)$ by setting $g\cdot\varphi=(g^{-1})^\ast (\varphi)$ 
for any $g\in\Gamma$ and $\varphi\in C^n (\xtil)$, where $g^\ast=\,^t\! g_\ast$ is the usual map induced
by $g$ on cochains.  
This action leaves $C_c^\ast (\xtil)$, ${C}_b^\ast (\xtil)$ and ${C}_{b,c}^\ast (\xtil)$ 
invariant and commutes with the differential, thus endowing these modules (and $C^\ast (\xtil)$, of course)
with a $\Gamma$--complex structure.

It is proved in~\cite{Ivanov}
that for every $n\geq 0$ the $\Gamma$--module ${C}_b^n (\xtil)$ is relatively injective.
We show here that the same is true for $\cclim^n (\xtil)$, and that the modules $\cc^n (\xtil)$,
$C^n (\xtil)$ are injective.

\begin{prop}\label{inj1:prop}
Let $n\geq 0$.
The $\Gamma$--modules $C^n (\xtil)$ and $\cc^n (\xtil)$ are injective.
The Banach $\Gamma$--modules $\clim^n (\xtil)$ and $\cclim^n (\xtil)$ are relatively injective.
\end{prop}
\noindent{\sc Proof:}
Let $\iota\colon A\to B$ be an injective map between unbounded $\Gamma$--modules,
with left inverse $\sigma\colon B\to A$, and suppose we are given a 
$\Gamma$--map $\alpha\colon A\to C^n (\xtil)$. 
We denote by $e_0,\ldots,e_n$ the vertices of the standard $n$--simplex, and define
$\beta\colon B\to C^n (\xtil)$ as follows: given $b\in B$, the cochain $\beta (b)$
is the unique linear extension of the map that on the singular
simplex $s$ takes the following value:
$$
\beta (b) (s) =\sum_{g\in \Gamma} h_{\widetilde{X}} \left(g^{-1} (s(e_0))\right)
\cdot \left(\alpha (g(\sigma (g^{-1}(b))))(s)\right),
$$
where $h_{\widetilde{X}}$ is the map provided by Lemma~\ref{tecnico}.
By Lemma~\ref{tecnico}--(1), the sum involved is in fact finite, so $\beta$ is well-defined.
Moreover, for every $b\in B$, $g_0 \in\Gamma$ and $s\in S_n (\xtil)$ we have
$$
\begin{array}{lll}
\beta (g_0\cdot b)(s) & = & \sum_{g\in \Gamma} \hx (g^{-1}(s(e_0)))\cdot (\alpha (g(\sigma(g^{-1}
g_0 (b))))(s))\\ &=&
\sum_{g\in \Gamma} \hx (g^{-1}g_0 (g_0^{-1}\cdot s)(e_0)))\cdot (\alpha (g_0 g_0^{-1}g(\sigma(g^{-1}
g_0 (b))))(s))\\ &=&
\sum_{k\in \Gamma} \hx (k^{-1} (g_0^{-1}\cdot s)(e_0)))\cdot (\alpha (g_0 k(\sigma(k^{-1}
(b))))(s))\\ &=&
\sum_{k\in \Gamma} \hx (k^{-1} (g_0^{-1}\cdot s)(e_0)))\cdot (\alpha (k (\sigma(k^{-1}
(b))))(g_0^{-1}\cdot s))\\ &=&
\beta (b) (g_0^{-1}\cdot s)= (g_0 \cdot \beta (b))(s),
\end{array} 
$$
so $\beta$ is a $\Gamma$--map. Finally, 
$$
\begin{array}{lll}
\beta(\iota(b))(s)&=& \sum_{g\in \Gamma} \hx (g^{-1}(s(e_0)))\cdot 
(\alpha (g(\sigma(g^{-1}(\iota(b)))))(s))\\ &=& 
\sum_{g\in \Gamma} \hx (g^{-1}(s(e_0)))\cdot 
(\alpha (g(\sigma(\iota(g^{-1}\cdot b))))(s))\\ &=&
\sum_{g\in \Gamma} \hx (g^{-1}(s(e_0)))\cdot 
(\alpha (b)(s))\\ &=&
\left(\sum_{g\in \Gamma} \hx (g^{-1}(s(e_0)))\right)\cdot (\alpha (b)(s))=\alpha(b)(s),
\end{array}
$$
so $\beta\circ\iota=\alpha$.
Thus $C^n (\xtil)$ is an injective $\Gamma$--module.

Let now $s\in S_n (\xtil)$ be any singular $n$--simplex. By Lemma~\ref{tecnico}--(1)
there exists a neighbourhood $U$ of $s$ in $S_n (\xtil)$ such that the set
$\{g\in\Gamma\, | \, h_{\widetilde{X}} (g^{-1} (s' (e_0)))\neq 0\ {\rm for\ some}\
s'\in U\}$ is finite. This readily implies that if $\alpha (A)\subseteq C_c^n (\xtil)$, 
then also $\beta(B)\subseteq C_c^n (\xtil)$. Thus also $C_c^n (\xtil)$ is an injective 
$\Gamma$--module.

The same argument applies \emph{verbatim} 
if $C^n (\xtil)$ is replaced by $\clim^n (\xtil)$, and $A,B$ are Banach modules: moreover, it is easily seen
that if $\alpha$ is bounded and $||\sigma||\leq 1$, then also $\beta$ is bounded, and
$||\beta||\leq ||\alpha||$. This, together with the argument above about continuity, implies
that $\clim^n (\xtil)$ and $\cclim^n (\xtil)$ are relatively injective Banach $\Gamma$--modules.
\qed\smallskip
 
\subsection{Singular cochains as strong resolutions of $\R$}\label{conv:subsec} 
If $E^\ast$ is one of the complexes $\cst (\xtil)$, $\ccst (\xtil)$, $\climst (\xtil)$,
$\cclimst (\xtil)$, endowed with the usual structures of $\Gamma$--complexes, then
a natural map $\delta^{-1}\colon E^{-1}:=\R\to E^0$ is defined, 
such that for any $t\in\R$ and $x_0\in S_0 (\xtil)$ we have $\delta^{-1}(t) (x_0)=t$.
It is readily seen that $\delta^0\circ \delta^{-1}=0$, so 
the augmented sequence of modules thus obtained, which will still be denoted
by $({E}^\ast,\delta^\ast)$ from now on, is
a complex. We would now like to show that in some cases such an augmented 
complex is exact (for example, this is obviously true 
if $E^\ast=C^\ast (\xtil)$ and $\xtil$ is contractible), and moreover admits a contracting homotopy.

\begin{prop}\label{strong1:prop}
Suppose $\xtil$ is contractible. Then
the complexes ${\cst (\xtil)}$ and ${\ccst (\xtil)}$ are strong resolutions of
$\R$ as an unbounded $\Gamma$--module. Moreover, the complexes ${\climst (\xtil)}$ and 
${\cclimst (\xtil)}$ are strong resolutions of
$\R$ as a Banach $\Gamma$--module.
\end{prop}
\noindent {\sc Proof:} 
Since $\xtil$ is contractible, there exist $x_0\in\xtil$ and  a continuous map
$H\colon \xtil\times [0,1]\to\xtil$ such that $H(x,0)=x$
and $H(x,1)=x_0$ for every $x\in\xtil$.
For $n\geq 0$,  
let $e^n_0,\ldots,e^n_n$ be 
the vertices of the standard simplex $\Delta_n\subset \R^{n+1}$, 
and let $Q^n_0$ be the face of $\Delta_n$ opposite to $e_0^n$. 
Let also $r_n\colon Q^{n+1}_0\to \Delta_n$ be defined by
$r_n (t_1 e^{n+1}_1+\ldots t_{n+1} e^{n+1}_{n+1})=t_1 e^n_0+\ldots t_{n+1} e^n_n$. 
For $n\geq 0$, we define $T_n\colon C_n (\xtil)\to C_{n+1}(\xtil)$
as the unique linear map such that 
if $s\in S_n(\xtil)$, then the following holds:
if
$p=t e^{n+1}_0+(1-t)q \in\Delta_{n+1}$, where $q\in Q^{n+1}_0$, then $(T_n(s)) (p)=H(s(r_n(q)), t)$.
($T_n (s)$ is just the ``cone'' over $s$ with vertex $x_0$, contructed 
by using the contracting homotopy $H$).
By Lemma~\ref{basi:top:lemma3}, $T_n (s)$ is well-defined and continuous. 
Moreover, we define $T_{-1}\colon \R\to C_0 (\xtil)$
by $T_{-1} (t)=t x_0$.
It is readily seen that, if $d_\ast$ is the usual (augmented)
differential on singular chains, then $d_{0} T_{-1}={\rm Id}_\R$, and
for every $n\geq 0$ we have
$T_{n-1} \circ d_n + d_{n+1} \circ T_n= {\rm Id}_{C_n (\xtil)}$.

For every $n\geq 0$, let now $k^n\colon C^n (\xtil)\to C^{n-1} (\xtil)$ be defined by $k^n (\varphi)(c)=
\varphi (T_{n-1} (c))$. It is readily seen that $\{k^n\}_{n\in\matN}$ provides a contracting
homotopy for the complex 
${C^\ast (\xtil)}$, which is therefore a strong resolution of $\R$.
By Lemma~\ref{basi:top:lemma3}, the map 
$T_n |_{S_n (\xtil)}\colon S_n (\xtil)\to S_{n+1} (\xtil)$ is continuous,
so the contracting homotopy $\{k^n\}_{n\in\matN}$
restricts to a contracting homotopy for the augmented complex of continuous cochains 
${\ccst (\xtil)}$, which therefore also gives a strong resolution of $\R$.
Moreover, since $T_n$ sends a simplex to a simplex, if $\alpha\in \clim^n (\xtil)$ 
then  $||k^n (\alpha)||\leq ||\alpha||$. Thus a suitable restriction of $k^\ast$ 
provide contracting homotopies for the complexes of Banach $\Gamma$--modules 
${\climst (\xtil)}$ and ${\cclimst (\xtil)}$. These complexes give therefore
strong resolutions of $\R$ as a Banach $\Gamma$--module.
\qed\smallskip

The following result is very deep, and plays a fundamental r\^ole in the study of bounded
cohomology of topological spaces. 
Together with a separate argument providing the required
control on seminorms, 
it implies, for example, that the bounded cohomology
of a countable CW--complex is canonically isomorphic to the bounded cohomology of its fundamental
group~\cite[Section 3.1]{Gromov}.

\begin{teo}[\cite{Ivanov}]\label{fund:teo}
Suppose $X$ has the homotopy type of a path connected countable CW--complex. Then $\climst (\xtil)$
is a relatively injective strong resolution of $\R$ as a Banach $\Gamma$--module.
\end{teo}

We now come to the following important:

\begin{prop}\label{GtoF:prop}
Let ${F}^\ast (\Gamma)$ (resp.~$F_b^\ast (\Gamma)$) be the standard resolution of $\R$
as an unbounded (resp.~Banach) $\Gamma$--module. 
There exists a chain map $\beta^\ast\colon {F}^\ast (\Gamma) \to \ccst (\xtil)$ 
which extends the identity of $\R$ and is such that $\beta^n$ is norm--decreasing for every $n\in\matN$.
In particular, $\beta^\ast$ restricts to a chain map $\beta_b^\ast\colon {F}_b^\ast (\Gamma) \to \cclimst (\xtil)$ 
which extends the identity of $\R$ and is such that $||\beta^n||\leq 1$ for every $n\in\matN$.
\end{prop}
\noindent {\sc Proof:} 
Let $n\geq 0$.
For every $f\in {F}^n (\Gamma)$, we define $\beta^n (f)$ to be the unique
singular cochain such that for every $s\in S_n (\xtil)$ we have
$$
\beta^n (f) (s)=\sum_{(g_0,\ldots,g_n)\in \Gamma^{n+1}} h_{\widetilde{X}} (g_0^{-1} (s(e_0)))\cdot \ldots\cdot
h_{\widetilde{X}} (g_n^{-1} (s(e_n)))\cdot f(g_0,\ldots, g_n),
$$
where $h_{\widetilde{X}}\colon \xtil\to \R$ is the continuous map provided by Lemma~\ref{tecnico}.
It is readily seen that the sum involved in the definition above is finite, so 
$\beta^n (f) (s)$ is well-defined. Moreover, Lemma~\ref{tecnico}--(1) ensures that
for every $s\in S_n (\xtil)$ 
a neighbourhood $U$ of $s$ exists such that
$$
\{(g_0,\ldots,g_n)\in \Gamma^{n+1}\, | \, \exists s'\in U\quad {\rm s.t.}\quad  \hx (g_i^{-1} (s' (e_i)))\neq 0\ 
\forall i=0,\ldots,n\}
$$ 
is finite. This easily implies that $\beta^n (f)$ is indeed continuous. 

The fact that $\beta^n$ is norm--decreasing for every $n\in\matN$ is immediate,
and it is straightforward to check that $\beta^\ast$
is indeed a $\Gamma$--equivariant chain map.
\qed\smallskip

\section{Proofs of Theorems~\ref{lim2:teo}, \ref{lim:teo}, \ref{lowdim2:teo}}\label{main:sec}
We are now ready to deduce Theorems~\ref{lim2:teo},~\ref{lim:teo}
from the results about resolutions proved in the preceding section.
Looking closely at the formula involved in the proof
that $C_c^n (\widetilde{X})$ is relatively injective, and 
at the explicit description for the contracting homotopy
for $C^\ast (\widetilde{X})$ when $\widetilde{X}$ is contractible, we will be able to write down the explicit
formulae required in the statement of Theorem~\ref{lim:teo}.
Moreover, in Theorem~\ref{lowdim2:teo}, 
similar formulae will be obtained  in the one-dimensional case
even without the assumption that $\widetilde{X}$ is contractible. 

We begin with the following: 

\begin{lemma}\label{sollevo:lemma}
For any topological space $X$,
the chain map
$p^\ast\colon C^\ast (X)\to C^\ast (\xtil)$ restricts to the following isometric isomorphisms
of complexes:
$$
\begin{array}{ll}
p^\ast\colon C^\ast (X)\to C^\ast (\xtil)^\Gamma,\quad & 
p^\ast|_{C_c^\ast(X)}\colon C_c^\ast (X)\to C_c^\ast (\xtil)^\Gamma, \\
p^\ast|_{{C}_b^\ast(X)}\colon \climst (X)\to \climst (\xtil)^\Gamma, \quad &
p^\ast|_{\cclimst (X)}\colon \cclimst (X)\to \cclimst (\xtil)^\Gamma,
\end{array}
$$
which induce therefore isometric isomorphisms
$$
\begin{array}{ll}
\hst (X) 
\cong H^\ast (\cst (\xtil)^\Gamma), \quad &  
\hcst (X)
\cong H^\ast (\ccst (\xtil)^\Gamma),\\
\hlimst (X)
\cong H^\ast (\climst (\xtil)^\Gamma), \quad &  
\hclimst (X)
\cong H^\ast (\cclimst (\xtil)^\Gamma).
\end{array}
$$
\end{lemma}
\noindent {\sc Proof:}
The fact that $p^\ast$ is an isometric embedding on the space of 
$\Gamma$--invariant cochains on $\xtil$ is obvious, thus the only non-trivial issue
to prove is the fact that $p^\ast (\varphi)$ is continuous if and only if $\varphi$
is continuous. 
By Lemma~\ref{open:lemma}, the map $p_\ast\colon S_n(\xtil)\to S_n (X)$ is a covering.
In particular, it is continuous, open and surjective, and this readily implies
that if $\varphi\colon S_n (X)\to\R$ is any map, then $\varphi$ is continuous if and only
if $\varphi\circ  p_\ast\colon S_n (\xtil)\to\R$ is continuous, whence the conclusion.
\qed\smallskip

\subsection{Proof of Theorem~\ref{lim2:teo}}
Suppose now that $X$ is a path connected countable CW--complex.
The inclusion $\widetilde{i}^\ast \colon {\ccst}(\xtil)\to{\cst}(\xtil)$
is a norm--decreasing chain $\Gamma$--map, and induces therefore
a norm--decreasing chain $\Gamma$--map
$\widetilde{i}^\ast_b \colon {\cclimst}(\xtil)\to{\climst}(\xtil)$.
We would like to show that
$\widetilde{i}^\ast_b$ admits a norm--decreasing right homotopy inverse.

By Theorem~\ref{fund:teo}, 
the complex $\climst(\xtil)$ provides a relatively injective strong resolution of
$\R$ as a Banach $\Gamma$--module. Thus Proposition~\ref{norm:prop} shows that there exists
a norm--decreasing chain map $\alpha_b^\ast\colon  \climst (\xtil)\to {F}_b^\ast(\Gamma)$ 
which extends the identity of $\R$. Let now $\beta_b^\ast \colon
{F}_b^\ast (\Gamma)\to \cclimst(\xtil)$ be the chain map provided by Proposition~\ref{GtoF:prop}.
By Propositions~\ref{inj1:prop} and~\ref{ext:prop} the map $\widetilde{i}^\ast_b\circ (\beta_b^\ast\circ\alpha_b^\ast)$
is $\Gamma$--homotopic to the identity of $\climst (\xtil)$. 
Since both $\widetilde{i}^\ast_b$ and 
$\beta_b^\ast\circ\alpha_b^\ast$ are norm--decreasing, the map 
$\widetilde{i}_b^\ast$ restricts to a map
$$
\widetilde{i}_b^\ast|_{\cclimst (\xtil)^\Gamma}\colon \cclimst (\xtil)^\Gamma\to \climst (\xtil)^\Gamma
$$
which induces in homology a norm--decreasing map $\overline{i}_b^\ast$
admitting a norm--decreasing right inverse (such an inverse is therefore an isometric embedding). 
Moreover, under the isometric identifications 
$\hlimst (X)\cong H^\ast (\climst (\xtil)^\Gamma)$,  
$\hclimst (X)\cong H^\ast (\cclimst (\xtil)^\Gamma)$ provided by 
Lemma~\ref{sollevo:lemma}, the map $\overline{i}_b^\ast$ corresponds to
$i_b^\ast \colon \hclimst (X)\to \hlimst (X)$,
whence the conclusion.
\qed\smallskip

\begin{rem}
If we were able to show that the complex $\cclimst (\xtil)$ provides a strong resolution
of $\R$, we could prove that, under the hypothesis of Theorem~\ref{lim2:teo}, 
the map $i^\ast_b$ is an isometric isomorphism.
However, it is not clear to us why the contracting homotopy for $\climst (\xtil)$ constructed in~\cite[Theorem 2.4]{Ivanov}
should take continuous cochains to continuous cochains, thus restricting to a contracting homotopy
for $\cclimst (\xtil)$. 
\end{rem}

\subsection{Proof of Theorem~\ref{lim:teo}}\label{standproof}
Let $\widetilde{\theta}^\ast\colon \cst(\xtil)\to \ccst (\xtil)$
be defined as the composition $\beta^\ast\circ\alpha^\ast$, where
$\alpha^\ast$, $\beta^\ast$
are the maps described in Remark~\ref{unbounded:rem}
and Proposition~\ref{GtoF:prop}.
Since $\xtil$ is contractible,
by Propositions~\ref{inj1:prop} and~\ref{strong1:prop}
both $\cst (\xtil)$ and $\ccst (\xtil)$
provide injective strong resolutions of
$\R$. Therefore, 
by Proposition~\ref{ext:prop} the compositions
$\widetilde{i}^\ast\circ \widetilde{\theta}^\ast$ and $\widetilde{\theta}^\ast\circ \widetilde{i}^\ast$ are $\Gamma$--homotopic to
the identity respectively of $\cst (\xtil)$ and of $\ccst (\xtil)$. 
Therefore, $\widetilde{i}^\ast$, $\widetilde{\theta}^\ast$
restrict to homotopy equivalences between $\ccst (\xtil)^\Gamma$ and $\cst (\xtil)^\Gamma$,
which in turn define isomorphisms 
$$
\overline{i}^\ast\colon
H^\ast (\ccst (\xtil)^\Gamma)\to H^\ast (\cst (\xtil)^\Gamma),\quad
\overline{\theta}^\ast\colon
H^\ast (\cst (\xtil)^\Gamma)\to H^\ast (\ccst (\xtil)^\Gamma),$$
that are one the inverse of the other.
Finally, under the identifications 
provided by 
Lemma~\ref{sollevo:lemma}, the map $\overline{i}^\ast$ corresponds to
$i^\ast \colon \hcst (X)\to \hst (X)$, while $\overline{\theta}^\ast$ corresponds
to the inverse $\theta^\ast\colon \hst(X)\to \hcst(X)$ of $i^\ast$.
This proves in particular that $i^\ast$ is an isomorphism.
Moreover, since $\alpha^\ast$, $\beta^\ast$ are norm--decreasing,
so is $\theta^\ast$, and this readily implies that
$i^\ast$ is isometric.

In order to write down an explicit formula for
$(i^\ast)^{-1}=\theta^\ast$ it is sufficient to exhibit an explicit construction of $\widetilde{\theta}$.
To this aim, let 
$H\colon \xtil\times [0,1]\to\xtil$ be a homotopy such that $H(x,0)=x$
and $H(x,1)=x_0$ for every $x\in\xtil$, and let $T_\ast\colon C_\ast (\xtil)\to
C_{\ast+1} (\xtil)$ be the induced contracting homotopy,
as described in the proof of Proposition~\ref{strong1:prop}.
Given $n\geq 0$ and $(g_0,g_1,\ldots,g_n)\in\Gamma^{n+1}$, we will now construct a simplex
$\overline{s}(g_0,g_1,\ldots,g_n)\in S_n (\xtil)$ such that
$\overline{s}(g_0,g_1,\ldots,g_n)(e_i)=g_i (x_0)$ 
for every $i=0,1,\ldots,n$. 
If $n=0$, we set $\overline{s}(g_0)=g_0 (x_0)$, while if
$\overline{s}(g_0,g_1,\ldots,g_n)\in S_n (\xtil)$ has been
defined for every $(g_0,g_1,\ldots,g_n)\in\Gamma^{n+1}$, then for every 
$(g_0,g_1,\ldots,g_{n+1})\in\Gamma^{n+2}$
we set 
$$
\overline{s}(g_0,g_1,\ldots,g_{n+1})=g_0\cdot(T_n (g_0^{-1}\cdot \overline{s}(g_1,g_2,\ldots,g_{n+1})))
$$
(so $\overline{s}(g_0,g_1,\ldots,g_{n+1})$ is the image under $g_0$ of the 
``cone'' of vertex $x_0$
based on $\overline{s}(g_0^{-1}g_1,\ldots,g_0^{-1}g_{n+1})$).

It is not difficult to check that, according to Propositions~\ref{strong1:prop},~\ref{GtoF:prop} and Remark~\ref{unbounded:rem}, the map
$\widetilde{\theta}^\ast=\beta^\ast\circ \alpha^\ast \colon \cst (X)\to \ccst (X)$
is defined as follows: 
if $\widetilde{\varphi}\in C^n (\xtil)$, then $\widetilde{\theta}^n(\widetilde\varphi)\in C^n_c (\xtil)$ is the unique cochain
such that for every $\widetilde{s}\in S_n (\xtil)$ we have 
$$
\widetilde{\theta}^n (\widetilde\varphi) (\widetilde{s})=
\sum_{(g_0,g_1,\ldots,g_n)\in\Gamma^{n+1}} \hx (g_0^{-1} (\widetilde{s}(e_0)))\cdot \ldots\cdot
\hx (g_n^{-1} (\widetilde{s}(e_n)))\cdot \widetilde{\varphi} (\overline{s}(g_0,\ldots,g_n)).
$$
As a consequence, 
for every $[\varphi]\in H^n (X)$,  the coclass $(i^n)^{-1} ([\varphi])=\theta^n ([\varphi])$
is represented by the cocycle that maps every $s\in S_n (X)$ to the real number
$\widetilde{\theta}^n (p^\ast (\varphi)) (\widetilde{s})$, where $\widetilde{s}$ 
is any lift of $s$ to $\xtil$.

Being norm--decreasing in every degree, the map $\widetilde{\theta}^\ast\colon \cst (\xtil)\to \ccst (\xtil)$ restricts to
a chain map $\widetilde{\theta}_b^\ast\colon \climst (\xtil)\to \cclimst(\xtil)$ extending the identity
of $\R$. By Propositions~\ref{strong1:prop} and~\ref{inj1:prop}, this map provides a
$\Gamma$--homotopy inverse to $\widetilde{i}^\ast_b$, and induces as above a norm--decreasing
inverse $\theta^\ast_b\colon \hlimst(X)\to\hclimst(X)$ of $i^\ast_b$.
This proves that $i^\ast_b$ is an isometric isomorphism.
Moreover, for every $[\varphi_b]\in H_b^n (X)$,  the coclass $(i_b^n)^{-1} ([\varphi_b])=\theta_b^n ([\varphi_b])$
is represented by the cocycle that maps every $s\in S_n (X)$ to the real number
$\widetilde{\theta}_b^n (p^\ast (\varphi)) (\widetilde{s})$, where $\widetilde{s}$ 
is any lift of $s$ to $\xtil$.


\subsection{Proof of Theorem~\ref{lowdim2:teo}}\label{lowdim2:sub}
In this subsection we concentrate on the relations between continuous
and ordinary cohomology in dimension one.
We begin by dealing with the map $i^1\colon H^1_c (X)\to H^1 (X)$ between \emph{unbounded}
cohomology modules.
Since $X$ is locally path connected, 
by Proposition~\ref{connected:prop} we can suppose that $X$ is path connected.

We first show that $i^1$ is surjective, \emph{i.e.}~that
any cocycle $\varphi\in C^1 (X)$ is cohomologous to a continuous cocycle. So, let 
$\widetilde{\varphi}=p^\ast (\varphi)\in C^1 (\xtil)$. Since $\varphi$ is a cocycle and $\xtil$
is simply connected, $\widetilde{\varphi}$ is a coboundary, so $\widetilde{\varphi}=\delta f$
for some $f\in C^0 (\xtil)$. Moreover, for every $g\in\Gamma$ we have $g\cdot\widetilde{\varphi}=\widetilde{\varphi}$,
and this readily implies that
\begin{equation}\label{equi:eq}
f(g(y)) -f(g(x))=f(y)-f(x)\quad {\rm for\ every}\ x,y\in\xtil,\ g\in\Gamma .
\end{equation}
 
We now replace $f$ with the map
$f_c\colon \xtil\to\R$
defined by 
$$
f_c (x)=\sum_{g\in\Gamma} h_{\widetilde{X}}
(g^{-1}(x))\cdot f(g(x_0)),
$$
where $x_0\in\xtil$ is a fixed basepoint,
and $h_{\widetilde{X}}$ is the function provided by Lemma~\ref{tecnico}.
By Lemma~\ref{tecnico}--(1), if $x\in \xtil$ there exists a neighbourhood
$U$ of $x$ in $\xtil$ such that the set $\{g\in \Gamma\, | \, h_{\widetilde{X}} (g^{-1}(x'))\neq 0
\ {\rm for\ some}\ x'\in U\}$ is finite. This readily implies that $f_c$ is well-defined and continuous,
and determines therefore a continuous $0$--cochain, which we will still denote by $f_c$.

Let us consider the difference $\widetilde k=(f-f_c)\colon\xtil\to\R$.
For every $g_0\in \Gamma$, $x\in\xtil$ we have
$$
\begin{array}{lll}
\widetilde{k}(g_0^{-1} (x))-\widetilde{k}(x) 
& = & f(g_0^{-1} (x))- f(x)-\sum_{g\in\Gamma} h_{\xtil} (g^{-1} (x))(f(gg_0^{-1} (x))- f(g (x)))\\
& = &
f(g_0^{-1} (x))- f(x)-\sum_{g\in\Gamma} h_{\xtil} (g^{-1} (x))(f(g_0^{-1} (x))- f(x))\\
& = &
f(g_0^{-1} (x))- f(x)-\left(\sum_{g\in\Gamma} h_{\xtil} (g^{-1} (x))\right) (f(g_0^{-1} (x))- f(x))\\
& = & 0,
\end{array}
$$
where the second equality is obtained by specializing
equation~\eqref{equi:eq} to the case $y=g_0^{-1} (x)$.
Being $\Gamma$--equivariant, the map $\widetilde k$ defines therefore
a unique map $k\colon X\to \R$ such that $\widetilde{k} (x)=k (p(x))$ 
for every $x\in \xtil$. We now set
$$
\varphi_c =\varphi -\delta k \in C^1 (X) .
$$
We have by construction $[\varphi_c ]=[\varphi]$ in $H^1 (X)$,
so in order to show that $i^1$ is surjective it is sufficient
to show that $\varphi_c$ is continuous, or, equivalently, that
$\widetilde{\varphi}_c=p^\ast (\varphi_c)\in C^1 (\xtil)$
is continuous (see Lemma~\ref{sollevo:lemma}). 
However, we have
$$
\widetilde{\varphi}_c=p^\ast (\varphi -\delta k)=
p^\ast (\varphi)-\delta p^\ast (k)= \widetilde{\varphi}
-\delta \widetilde{k}=\delta f- \delta (f-f_c)=\delta f_c ,
$$
which is continuous since $f_c\in C^0_c (\xtil)$.
We have thus proved that $i^1$ is surjective,
also providing a somewhat explicit procedure for 
replacing a singular $1$--cocycle with a cohomologous continuous
$1$--cocycle.

Let now $\varphi\in C^1_c (X)$ be a continuous cocycle with $i^1 ([\varphi])=0$, fix a basepoint
$x_0\in X$, and, for every $q\in X$, let $s_q\colon [0,1]\to X$
be a fixed continuous path with $s_q (0)=x_0$, $s_q (1)=q$. We define a real function
$f\colon X\to\R$ by setting $f(q)=\varphi (s_q)$. This function defines a $0$--cochain
which will still be denoted by $f$. We will now show that $\delta (f)=\varphi$, and that
$f$ is continuous, thus proving that $\varphi$ is the coboundary of a continuous
$0$--cochain, and that $i^1$ is injective. 

So, let $s\in S_1 (X)$ be a simplex with endpoints $q_0,q_1$. Since $s+s_{q_0}-s_{q_1}$
is a cycle and $\varphi$ is a coboundary, we have 
$$
\varphi (s)=\varphi (s_{q_1})-\varphi (s_{q_0})=f(q_1)-f (q_0)=\delta (f) (s),
$$
so $\varphi=\delta (f)$. We now show that $f$ is continuous.
Let $q\in X$ and $\varepsilon>0$ be given. If $c_q$ is the constant
$1$--simplex with $c_q (t)=q$ for every $t\in [0,1]$, then $\varphi (c_q)=f(q)-f (q)=0$. Thus
it is not difficult to show that since $\varphi$ is continuous
there exists a path connected neighbourhood
$U$ of $q$ such that $|\varphi (s)|< \varepsilon$ for every simplex $s$ with values in $U$.
In particular, for every $r\in U$ there exists a simplex $s_{q,r}$ such that $s_{q,r} (0)=q$,
$s_{q,r} (1)=r$, and $|\varphi (s_{q,r})|<\varepsilon$, so 
$|f(r)-f(q)|=
|\varphi (s_{q,r})|<\varepsilon$. We have thus shown that $f$ is continuous, so $i^1$
is injective. Also note that if $\varphi$ is supposed to be bounded, then $||f||\leq ||\varphi||$,
so $i_b^1\colon {H}^1_{b,c} (X)\to {H}_b^1 (X)$ is injective too.
Moreover, 
the last statement is obviously still true even when $X$ is not path connected (but still locally path
connected). 

The fact that ${H}_b^1 (X)=0$ is well-known and easy: a bounded $1$--cocycle $\omega$ defines a bounded
homomorphism between $\pi_1 (X)$ and $\R$, but $\R$ does not contain 
non-trivial bounded subgroups, so $\omega$
has to be the coboundary of a $0$--cochain $\psi$; moreover, 
the same argument showing above that $|| f ||\leq ||\varphi ||$ applies
here ensuring 
that $\psi$
can be chosen in such a way that $||\psi||\leq ||\omega||$, and this implies
that $[\omega]=0$ in ${H}_b^1 (X)$, so ${H}_b^1 (X)=0$.
Therefore also ${H}_{b,c}^1 (X)=0$, since 
$i_b^1\colon {H}^1_{b,c} (X)\to {H}_b^1 (X)$ is injective as noted above. 
\qed\smallskip

\section{Two (counter)examples}\label{1dim:sec}

It is not difficult to construct disconnected spaces whose continuous
cohomology is not isomorphic to standard cohomology, even in dimension
$0$.
For example,
if $Y$ is the Cantor set, then any real function on $Y$ defines a $0$--cocycle, 
while a $0$--cocycle is continuous if and only if it is defined by a continuous
real function on $Y$,
and this readily shows that
in this case the map $i^0\colon H^0_c (Y)\to H^0 (Y)$ is not
surjective.

The following examples provide path connected spaces whose continuous 
cohomology is 
not isomorphic (through $i^\ast$) to standard cohomology, even in dimension one.

\begin{center}
\begin{figure}
\input{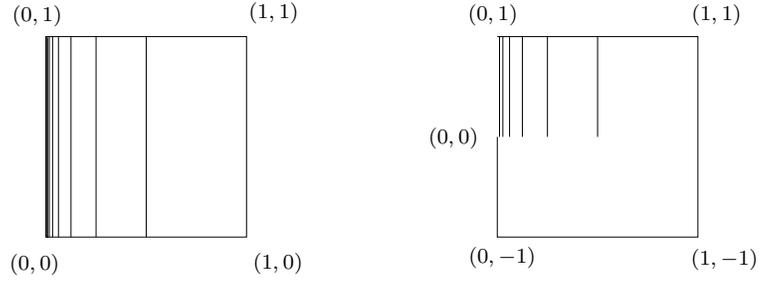}
\caption{The spaces $X_1$ (on the left) and $X_2$ (on the right).}\label{contro:fig}
\end{figure}
\end{center}

Let 
$$
X_1=\left( \left(\{0\} \cup \left(\bigcup_{n\geq 1}
\left\{\frac {1}{n}\right\}\right)\right)\times [0,1]\right) \cup \left([0,1]\times \{0,1\}\right)\subset \R^2
$$
be endowed with the Euclidean topology (see Figure~\ref{contro:fig}--left). 
For $i\in\matN$, let $\alpha_i\colon [0,1]\to X$ be the constant--speed
parameterization of the polygonal path with vertices $(0,0)$, $(1/i,0)$, $(1/i,1)$ and $(0,1)$
if $i\neq 0$,
and $\alpha_i (t)=(0,1-t)$ if $i=0$. It is not difficult to prove that
$H_1 (X_1)$
is freely generated (as a vector space) by the classes represented
by the loops $\alpha_i\ast \alpha_0$, $i\geq 1$. Thus, if $i^1$ were surjective,
by the Universal Coefficient Theorem,
for any real sequence $\{\epsilon_i\}_{i\geq 1}$
there should exist a continuous cocycle $\psi\in C_c^1 (X_1)$ such that 
$\psi (\alpha_0)+\psi (\alpha_i)=\epsilon_i$.
Now choosing $\epsilon_i= (-1)^i$ we would get 
$$
|\psi (\alpha_i)-\psi (\alpha_{i+1})|= |\psi (\alpha_i)+\psi (\alpha_0)-(\psi (\alpha_{i+1})+
\psi (\alpha_0))|=2,
$$
and this would contradict the continuity of $\psi$, since $\lim_{i\to \infty} \alpha_i =\alpha_0^{-1}$
in the compact--open topology, where $\alpha_0^{-1}\colon [0,1]\to X_1$
is defined by $\alpha_0^{-1} (t)=\alpha_0 (1-t)$
for every $t\in [0,1]$.

It is maybe worth mentioning that as a byproduct of our results we obtain that $X_1$ 
does not admit a universal covering:
in fact, $X_1$ is metrizable, whence paracompact, while the surjectivity of $i^1$ was 
established in the proof of Theorem~\ref{lowdim2:teo} without using any local path connectedness
assumption.

Let now $X_2\subset \R^2$ be defined as follows (see Figure~\ref{contro:fig}--right):
$$
\begin{array}{ccccc}
X_2 &=& \left( \left(\bigcup_{n\geq 1}
\left\{\frac {1}{n}\right\}\right)\times [0,1]\right) &\cup&  \left( (0,1]\times \{-1,1\}\right)\\
&\cup & (\{0\}\times [-1,0])& \cup & (\{1\}\times [-1,1]).
\end{array}
$$
We let $f\colon X_2\to \R$ be the $0$--cochain such that for every $(x,y)\in X_2$ 
$$
f(x,y)=\left\{
\begin{array}{cl}
0\quad & {\rm if}\   x=0,\, x=1,\ {\rm or}\ y=-1 \\
1-y \quad & {\rm otherwise}   
\end{array}
\right. 
$$
and set $\varphi=\delta f\in C^1 (X_2)$. The map $f$ is continuous at every point of $X_2$
except $(0,0)$, and we now sketch a proof of the fact that
$\varphi$ is indeed continuous. Since $X_2$ is metrizable and $[0,1]$ is compact,
it is sufficient to show that $\varphi$ is sequentially continuous. So, let us suppose by contradiction
that 
$\lim_{n\to\infty} s_n = s\in S_1 (X_2)$, 
while $\varphi (s_n)=f(s_n (1))-f(s_n(0))$ does not tend to $\varphi (s)=f(s(1))-f(s(0))$.
Then $\llim f(s_n (t_0))\neq f(s(t_0))$ for some $t_0\in\{0,1\}$. Since $f$ is everywhere continuous 
 except that at $(0,0)$ we 
must have $s(t_0)=(0,0)$, and $x(s_n(t_0))>0$  for an infinite number of indices. But it is easily seen
that a sequence of 
continuous
paths starting (resp.~ending) in points of $X_2$ with positive $x$--coordinate can converge to a path starting
(resp.~ending)
in $(0,0)$ if and only if it converges to the constant path. However, if $s$ is the constant path
in $(0,0)$  
and $\llim s_n=s$, then it is easily seen that definitively $x(s_n (0))=x(s_n (1))$, so 
$\llim \varphi (s_n)=0=\varphi (s)$, a contradiction.

Thus $\varphi\in C^1_c (X_2)$, and $\delta\varphi=\delta^2 f=0$. We now claim that
$[\varphi]\neq 0$ in $H^1_c (X_2)$, thus showing that $i^1\colon H^1_c (X_2)\to
H^1 (X_2)$ is not injective, since $[\varphi]=[\delta f]=0$ in $H^1 (X_2)$. 
In fact, if $\varphi$ were the coboundary of a \emph{continuous} $0$--cochain
$f_c$, then we would have $\delta (f-f_c)=0$. Since $X_2$ is path connected,
this would imply that $k\in\R$ should exist such that $f(x)=f_c (x)+k$ for 
every $x\in X_2$, a contradiction since $f_c$ is continuous while $f$ is not. 

Note that this example does not contradict Theorem~\ref{lim:teo}. In fact,
even if it is simply connected and $1$--dimensional, $X_2$ is not contractible,
since by Proposition~\ref{hom:prop} any contractible space has trivial first continuous
cohomology group.

\section{Smooth cohomology}\label{smooth:sec}
Suppose now $X$ is a smooth manifold. 
As mentioned in the introduction, we would like to concentrate our attention only on \emph{smooth} simplices,
and on cochains which are continuous with respect to the $C^1$--topology on smooth simplices.
More precisely, we set ${_s S_q} (X)=\{s\in S_q (X)\, |\, s\ {\rm is\ smooth}\}$ and we
denote by ${_s C_\ast} (X)$ the subcomplex of $C_\ast (X)$ generated by $_s S_\ast (X)$;
moreover, we let ${_s C^q} (X)$ be the dual space of ${_s C_q} (X)$. Of course,
${_s C^\ast} (X)$ is a differential complex, whose elements will be called
\emph{smooth cochains}. The homology of ${_s C^\ast} (X)$ will be called 
\emph{smooth cohomology} of $X$ and denoted
by ${_s \hst} (X)$. We will also denote by ${_s \climst} (X)\subseteq {_s \cst} (X)$ 
the subcomplex 
of bounded cochains, and by ${_s \hlimst} (X)$ the corresponding cohomology.

The inclusion $j_\ast \colon {_s C_\ast} (X)\hookrightarrow C_\ast (X)$ induces restrictions
$j^\ast\colon \cst (X)\to {_s \cst (X)}$, ${j}_b^\ast\colon\climst (X)\to {_s\climst (X)}$, 
which induce in turn maps $r^\ast\colon \hst (X)\to {_s \hst (X)}$,
${r}^\ast_b\colon \hlimst (X)\to {_s \hlimst (X)}$. The following well-known result
shows that smooth cohomology is isometric to the usual one:

\begin{prop}\label{smoothing:prop}
The maps $r^\ast\colon \hst (X)\to {_s \hst (X)}$,
${r}_b^\ast\colon \hlimst (X)\to {_s \hlimst (X)}$ are isometric isomorphisms.
\end{prop}
\noindent {\sc Proof:}
It is well-known (see~\emph{e.g.}~\cite[Theorem 16.6]{Lee})
that there exist maps $l_\ast \colon C_\ast (X)\to {_s C_\ast (X)}$,
$T_\ast\colon C_\ast (X)\to C_{\ast+1} (X)$ 
with the following properties: for every $s\in S_n (X)$, $l_n (s)$
is a simplex in $_s S_n (X)$, while $T_n (s)$ is the sum (with signs)
of a fixed number
(depending only on $n$) of simplices in $S_{n+1} (X)$; 
$l_\ast\circ j_\ast={\rm Id}_{{_s C_\ast} (X)}$; 
$j_\ast\circ l_\ast-{\rm Id}_{{C_\ast}(X)}=d\circ T_\ast +T_\ast\circ d$.
The dual maps $l^\ast$, $j^\ast$ are therefore norm--decreasing, and one the $\Gamma$--homotopy inverse
of the other. Since $T_n$ is bounded for every $n$, the same is true for the restrictions
of $l^\ast$, $j^\ast$ to bounded cochains, whence the conclusion.
\qed\smallskip

For every $n\in\matN$, we now endow ${_s S_n} (X)$ with the $C^1$--topology 
(see Appendix~\ref{basi:app} for the definition and the needed properties 
of $C^1$--topology). We say that a cochain $\varphi\in {_s C^n (X)}$ is 
\emph{continuous} if it restricts
to a continuous map on ${_s S_n (X)}$, and we denote by ${_s \ccst (X)}$
the subcomplex of continuous cochains in ${_s \cst} (X)$.
We also denote by ${_s \cclimst} (X) = {_s \ccst (X)}\cap {_s \climst (X)}$
the complex of bounded continuous cochains.
The corresponding cohomology modules will be denoted by
${_s \hcst (X)}$ and ${_s \hclimst (X)}$.
The natural inclusions of cochains induce maps
$$
{_s i^\ast} \colon  {_s \hcst (X)}\to {_s \hst (X)},
\qquad 
{_s i_b^\ast}  \colon  {_s \hclimst (X)}\to {_s \hlimst (X)}.
$$

Basically, all the results proved in the preceding sections for continuous and
singular cohomology of sufficiently nice topological space extend to the cohomology theories
just introduced for smooth manifolds. We state here the 
facts we will need in 
Section~\ref{gromov:sec}, also giving an outline of their proofs.

\begin{lemma}\label{smoothinv:lemma}
The map ${_s i_b^\ast}$ admits a norm--decreasing right inverse.
\end{lemma}
\noindent
{\sc Proof:}
Since the $C^1$--topology
is finer than the compact--open topology, the restriction map $j^\ast\colon C^\ast (X)\to {_s C^\ast (X)}$
introduced above
takes continuous cochains into continuous smooth cochains.
Therefore $j^\ast$ restricts to a chain map
$\cclimst (X)\to {_s \cclimst (X)}$, which induces in turn
a map ${r}_{b,c}^\ast\colon \hclimst (X)\to {_s \hclimst (X)}$. We have therefore the sequence of maps
$$
\xymatrix{
\hclimst (X) \ar[r]^{{r}_{b,c}^\ast} 
& {_s \hclimst (X)} \ar[r]^{_s i_b^\ast} & {_s \hlimst (X)} \ar[r]^{({r}_b^\ast)^{-1}} &
\hlimst (X) \ar@/^1.2pc/[lll]^{\beta^\ast}
}
$$
where $\beta^\ast$ is the right inverse of $i_b^\ast$ provided by Theorem~\ref{lim2:teo}.
By the very definitions we have ${_s i_b^\ast} \circ r_{b,c}^\ast=
{r}_b^\ast\circ {i_b^\ast}$, so 
$$
\begin{array}{llllll}
& {_s i_b^\ast}\circ \left({{r}_{b,c}^\ast}\circ \beta^\ast
\circ {({r}_b^\ast)^{-1}}\right)&=& 
\left({_s i_b^\ast}\circ {{r}_{b,c}^\ast}\right)\circ \beta^\ast
\circ {({r}_b^\ast)^{-1}} & = &
\left({r}_b^\ast\circ {i_b^\ast}\right)\circ\beta^\ast\circ
{({r}_b^\ast)^{-1}}\\
= & {r}_b^\ast\circ \left({i_b^\ast}\circ\beta^\ast\right)
\circ {({r}_b^\ast)^{-1}} & = &
{\rm Id}_{\hlimst (X)}. & &
\end{array}
$$
Therefore ${_s i_b^\ast}$ admits a norm--decreasing right
inverse.
\qed\smallskip

\begin{teo}\label{bors:teo}
The map ${_s i^\ast}\colon {_s \hcst (X)}\to {_s \hst (X)}$ is an isometric isomorphism.
\end{teo}
\noindent {\sc Proof:} 
All the arguments developed in Section~\ref{bor:sec} in order to prove
Theorem~\ref{stand:teo} apply \emph{verbatim} to smooth cohomology,
thus proving that ${_s i^\ast}$ is an isomorphism.
More precisely, since every smooth manifold
is locally \emph{smoothly} contractible, by Lemma~\ref{invhomsmooth} 
the graded \emph{sheaves} of smooth cochains and continuous
smooth cochains both provide resolutions of the constant sheaf $\widetilde{\R}$ on $X$.
Moreover, such sheaves admit a structure of modules over the sheaf of continuous functions,
and are therefore acyclic. Therefore, sheafified smooth cohomology is canonically isomorphic
to sheafified continuous smooth cohomology.
Now, in order to prove that singular cohomology is isomorphic to sheafified singular
cohomology, in Section~\ref{bor:sec} we only needed the following facts:
\begin{enumerate}
\item
the existence of a barycentric subdivision (co)operator taking continuous
locally zero 
cochains to continuous locally zero cochains (see Proposition~\ref{vannobene:prop});
\item
the fact that a locally defined continuous cochain could be extended to a global one
(see Lemma~\ref{borelext:lemma}).
\end{enumerate}
The proof of fact~(2) given in Section~\ref{first:sec} applies \emph{verbatim}
when restricting our attention to smooth simplices, endowed with the $C^1$--topology. Moreover, 
using Lemma~\ref{smoothbasi:lemma1} and the fact that the $C^1$--topology is finer than the compact--open
topology, it is readily seen that the barycentric subdivision (co)operators defined 
in Section~\ref{first:sec} take continuous
locally zero smooth cochains to continuous locally zero 
smooth cochains, so~(1) also holds. 

Since ${_s i^\ast}$ is obviously norm--decreasing, 
the fact that ${_s i^\ast}$ is an isometry is now a consequence of 
Lemma~\ref{smoothinv:lemma}.
\qed\smallskip


Let now $p\colon \xtil\to X$ be the smooth universal convering
of $X$
(see the Appendix for the definition and basic properties of smooth coverings). 
By Lemma~\ref{smoothbasi:lemma1}, the covering $p$ induces a well-defined map
$p^\ast\colon {_s \ccst} (X)\to {_s \ccst (\xtil)}$.   
Moreover, if $\Gamma\cong \pi_1 (X)$ is the group of the
covering automomorphisms of $p$, then
$\Gamma$ acts on $\xtil$ as a group of diffeomorphisms.
Therefore, as noted in the Appendix, $\Gamma$ also acts on ${_s C_c^\ast (\xtil)}$, and we will
denote by ${_s C_c^\ast} (\xtil)^\Gamma\subseteq {_s C_c^\ast (\xtil)}$ the subcomplex
of $\Gamma$--invariant continuous smooth cochains.

The following result easily follows from Lemma~\ref{opensmooth:lemma} (see also
the proof of Lemma~\ref{sollevo:lemma}):

\begin{lemma}\label{sollevosmooth:lemma}
The chain map
$p^\ast\colon {_s \ccst} (X)\to {_s \ccst} (\xtil)$ induces the isometric isomorphism
of complexes
$$
p^\ast|_{{_s \ccst} (X)}\colon {_s \ccst}  (X)\to {_s \ccst} (\xtil)^\Gamma,
$$
which induces in turn an isometric isomorphism ${_s \hcst (X)}\cong H^\ast ({_s \ccst (\xtil)^\Gamma})$.
\end{lemma}

\section{Gromov's proportionality principle}\label{gromov:sec}
Before going into the proof of the proportionality principle, we briefly describe
Gromov's original approach to the issue.

As mentioned in Subsection~\ref{simpl:sub}, bounded cohomology provides the natural ``dual''
theory to $L^1$--homology, and is therefore deeply related to the simplicial volume. 
More precisely, it is not difficult to show that if $X$ is a compact connected Riemannian
manifold, then ${\rm Vol} (X) / ||X|| $ equals the seminorm of the coclass of $H^n (X)$ represented
by the Riemannian volume form of $X$ (see Corollary~\ref{obvious:cor} below).
Keeping notations from the preceding section,
if $\Gamma \cong \pi_1 (X)$ is the group of the covering automorphisms 
of $\xtil$,  then the volume form of $X$ lifts to the volume form of $\xtil$, which is
$G$--invariant, where $G$ is the group of orientation--preserving isometries of $\xtil$.
Moreover, the seminorm of the volume form of $X$ is equal to the seminorm 
of the volume form of $\xtil$ in the homology of the appropriate complex
of $\Gamma$--invariant cochains.
An averaging process on cochains now allows to show 
that this seminorm is equal to the seminorm of
the volume form of $\xtil$ in the homology of $G$--invariant cochains.  
As a consequence, the r\^ole of $\Gamma$ turns out to be immaterial, and
 ${\rm vol} (X) / ||X||$ only depends on the geometry of $\xtil$. 

The argument just outlined is basically
Gromov's original approach to the proportionality principle~\cite[Section 2.3]{Gromov}.
However, as Gromov himself points out, in order to formally define the above mentioned
averaging process, one should restrict only to continuous (or at least Borel measurable)  cochains.
The fact that this assumption is harmless is a consequence of the (isometric) isomorphism
Theorem~\ref{standiso:teo}. 

Our exposition closely follows (and is in fact inspired by)
Bucher-Karlsson's argument (see~\cite{Bucher}). However, since the volume form is \emph{not}
continuous with respect to the compact--open topology (see Remark~\ref{noncont:rem} below), we work here
in the slightly different setting of continuous smooth cohomology,
endowing the space of smooth simplices with the $C^1$--topology, rather than with the compact--open topology.

\subsection{The duality principle}
From now on, we denote by $X$ a $n$--dimensional 
compact connected oriented Riemannian manifold with real fundamental class
$[X]_\R\in H_n (X)$. Recall that there exists a well-defined product
$\langle\, , \, \rangle\colon H^n (X)\times H_n (X)\to\R$, called \emph{Kronecker product},
such that $\langle [\varphi], [z]\rangle =\varphi (z)$ for every cocycle $\varphi\in C^n (X)$
and every cycle $z\in C_n (X)$.
We denote by $[X]^\R\in H^n (\R)\cong \R$ the \emph{fundamental coclass} of $X$, \emph{i.e.}~the
unique coclass in $H^n (X)$ such that $\langle [X]^\R,[X]_\R\rangle =1$. 
Also recall that we denote by $c^n\colon {H}_b^n (X)\to H^n (X)$ the comparison map induced
by the inclusion of cochains.
The following result, due to Gromov~\cite{Gromov}, is based on
Hahn-Banach Theorem, and is proved \emph{e.g.}~in~\cite[Proposition F.2.2]{BePe}:

\begin{prop}\label{duality:prop}
Let $||X||=||[X_\R]||$ be the simplicial volume of $X$. Then
$$
||X||=\frac{1}{||[X]^\R||}=\sup \left\{\frac{1}{||\varphi||},\, \varphi\in {H}_b^n (X),\, 
c^n(\varphi)=[X]^\R\right\}.
$$ 
In particular, $||X||=0$ if and only if $[X]^\R\notin {\rm Im}\, c^n$, \emph{i.e.}~if 
and only if ${\rm Im}\, c^n=0$. 
\end{prop}

\subsection{The volume form}
We define a 
map ${\rm Vol}_X\colon {_s S_n} (X)\to \R$ by setting
$$
{\rm Vol}_X (s)=\int_s \omega_X,
$$
where $\omega_X\in\Omega^n (X)$ is the 
volume form of $X$. Of course, ${\rm Vol}_X$ is continuous with respect to the $C^1$--topology
on ${_s S_n} (X)$, so its linear extension to smooth $n$--chains, which will still be denoted by ${\rm Vol}_X$,
defines an element in ${_s C^n_c} (X)$.
By Stokes' Theorem, such an element is closed, and defines therefore elements
$[{\rm Vol}_X]\in {_s H^n} (X)$,
$[{\rm Vol}_X]_c\in {_s H^n_c} (X)$.  

\begin{rem}\label{noncont:rem}
The following example shows that the cochain ${\rm Vol}_X$ is in general
\emph{not} continuous with respect to the compact--open topology. 
In fact, let us consider the Euclidean plane $\R^2$, endowed with the
usual volume form $dx_1 \wedge dx_2$. 
Let $Y=\R^2\setminus \{(0,0)\}$, and for 
for every $n\geq 1$ let 
$f_n \colon Y \to Y$ be the map which corresponds
to $z \mapsto z^n / (n |z|^{n-1})$ under the identification $Y \cong \matC\setminus \{0\}$
(this map is the composition of the rescaling of ratio $1/n$ with a map 
that ``wraps'' $Y$ around the origin $n$ times). An easy computation shows that
$f_n$ is an area-preserving local diffeomorphism of $Y$ onto itself, so if $s\in {_s S_2 (\R^2)}$ is any smooth 
simplex with ${\rm Im} (s)\subseteq Y$ and ${\rm Vol}_{\R^2} (s)=\alpha \neq 0$, then $s_n=f_n\circ s$  is a smooth simplex
such that ${\rm Vol}_{\R^2} (s_n)=\alpha\neq 0$. On the other hand, if ${\rm Im} (s)$ is contained
in the ball $B(0, R)\subseteq \R^2$ of radius $R$ centered at the origin, then ${\rm Im} (s_n)$ is contained in the ball
$B(0,R/n)$. This readily implies that $\lim_{n\to\infty} s_n=s_0$ in the compact--open topology,
where $s_0$ is the constant simplex with values in $\{0\}\subseteq \R^2$. 
Since ${\rm Vol}_{\R^2} (s_0)=0$, this shows that ${\rm Vol}_{\R^2}$ is not continuous
if we endow ${_s S_2 (\R^2)}$ with the compact--open topology.
\end{rem}

Let $r^n\colon H^n (X)\to {_s H^n (X)}$ be the map introduced
at the beginning of Section~\ref{smooth:sec}.

\begin{lemma}\label{contvol:lemma}
We have $(r^n)^{-1} ([{\rm Vol}_X])={\rm Vol} (X)\cdot [X]^\R$.
\end{lemma}
\noindent {\sc Proof:} Since $H^n (X)\cong \R$, we have 
$(r^n)^{-1} ([{\rm Vol}_X])=\langle (r^n)^{-1} ([{\rm Vol}_X]),
[X]_\R\rangle \cdot [X]^\R$. Moreover, it is well-known that 
the fundamental class of $X$ can be represented by the sum of the simplices in a
positively oriented smooth
triangulation of $X$. Evaluating the cohomology class $(r^n)^{-1} ([{\rm Vol}_X])$
on such a sum we get the sum of the volumes of the simplices of the triangulation, \emph{i.e.}~the volume
of $X$.
\qed\smallskip

\begin{cor}\label{obvious:cor}
We have
$$
\frac{||X||}{{\rm Vol} (X)}= \frac{1}{||[{\rm Vol_X}]_c||]}.
$$
\end{cor}
\noindent {\sc Proof:}
Since ${_s i^n} ([{\rm Vol}_X]_c)=[{\rm Vol}_X]$, 
by Proposition~\ref{smoothing:prop}, Theorem~\ref{bors:teo} and Lemma~\ref{contvol:lemma} we have
$$||[{\rm Vol_X}]_c||=||(r^n)^{-1} ({_s i^n} ([{\rm Vol}_X]_c))||={\rm Vol} (X)\cdot ||[X]^\R||.$$
Therefore, by Proposition~\ref{duality:prop} we get
$||X||=1/||[X]^\R||={\rm Vol}(X)/{||[{\rm Vol_X}]_c||]}$.
\qed\smallskip

From now on, we denote by $\xtil$ the Riemannian universal covering of $X$. 
By Corollary~\ref{obvious:cor}, the proportionality principle can be restated as follows:

\begin{teo}\label{newprop:teo}
The value $||[{\rm Vol_X}]_c||$ only depends on the isometry type of $\xtil$.
\end{teo}

Thus our efforts will be henceforth devoted to proving Theorem~\ref{newprop:teo}.

\subsection{The transfer map}
From now on, we denote
by $G$ the group of orientation--preserving isometries
of $\xtil$, and by $\Gamma\cong \pi_1 (X)<G$ the group of covering
automorphisms of $\xtil$. It is well-known
that $G$ admits a Lie group structure inducing the compact--open topology. 
Moreover, there exists on $G$ a left-invariant regular Borel measure $\mu_G$, which 
is called \emph{Haar measure} of $G$ and  
is unique up to scalars.
Since $G$ contains a cocompact subgroup, 
its Haar measure is in fact
also right-invariant~\cite[Lemma 2.32]{Sauer}.
Since $X\cong \xtil /\Gamma$ is compact, there exists a Borel subset $F\subseteq G$
with the following properties: $F$ contains exactly one representative for each class
in $\Gamma\backslash G$ and $F$ is relatively compact in $G$. We will call
such an $F$ a \emph{fundamental region} for $\Gamma$ in $G$.
From now on, we normalize the Haar measure
$\mu_G$ in such a way that $\mu_G (F)=1$.

In order to avoid too heavy notations, if $H$ is a subgroup of $G$
we set
${_s \hcst} (\xtil)^H= H^\ast ({_s \ccst} (\xtil)^H)$.
We also endow 
${_s \hcst}(\xtil)^H$ with the seminorm induced by 
${_s \ccst} (\xtil)^H$.
Recall that by Lemma~\ref{sollevosmooth:lemma} we 
have an isometric isomorphism 
${_s \hcst} (\xtil)^\Gamma\cong {_s \hcst} (X)$. 
The chain inclusion ${_s \ccst} (\xtil)^G\hookrightarrow 
{_s \ccst} (\xtil)^\Gamma$ induces a norm--decreasing map
$$
{\rm res}^\ast \colon  {_s \hbst}(\xtil)^G \longrightarrow {_s \hbst}(\xtil)^\Gamma\cong {_s \hbst} (X).
$$

Following~\cite{Bucher}, we will now construct a norm--decreasing left inverse 
of $\res$. We begin with the following:

\begin{lemma}\label{contact:lemma}
Let $s_0\in {_s S_\ast (\xtil)}$ be fixed. Then the map
$ 
G
\to {_s S_\ast} (\xtil)$ defined by 
$g\mapsto g\cdot s_0=g\circ s_0
$
is continuous.
\end{lemma}
\noindent {\sc Proof:}
Let us consider $G$ as a subset of the space $F_s (\xtil,\xtil)$ of smooth functions
from $\xtil$ to itself. 
Since the elements of $G$ are isometries,
the compact--open and the $C^1$--topology coincide on $G$
(see \emph{e.g.}~\cite[Theorem 5.12]{Loh0}). Therefore the conclusion follows from
Lemma~\ref{smoothbasi:lemma1}.
\qed\smallskip

Take now $\varphi\in {_s C_c^i} (\xtil)$ and $s\in {_s C_i} (\xtil)$, 
and consider the function $f_\varphi^s \colon G\to \R$
defined by $f_\varphi^s (g)=\varphi (g\cdot s)$. By Lemma~\ref{contact:lemma}, $f_\varphi^s$ is
continuous, whence bounded on the relatively compact subset $F\subseteq G$.
Therefore a well-defined
cochain $\tr^i (\varphi)\in {_s C^i} (\xtil)$ exists such that for every $s\in {_s S_i} (\xtil)$ we have
$$
\tr^i (\varphi) (s)=\int_F f_\varphi^s (g)\, d\mu_G (g) = \int_F \varphi (g\cdot s)\, d\mu_G (g) .
$$

\begin{prop}\label{funziona:prop}
The cochain $\tr^i (\varphi)$ is continuous. Moreover,
if $\varphi$ is $\Gamma$--invariant, then $\tr^i (\varphi)$ is $G$--invariant,
while if $\varphi$ is $G$--invariant, then $\tr^i (\varphi)=\varphi$.
\end{prop}
\noindent {\sc Proof:}
Let us define a 
distance $d_S$ on ${_s S_i} (\xtil)$ as follows. It is well-known that the Riemannian structure on $\xtil$
induces a distance $d_{T\xtil}$ on the tangent bundle $T\xtil$. Moreover, 
$d_{T\xtil}$
is $G$--invariant, in the sense that for every $g\in G$ the differential $dg\colon T\xtil\to T\xtil$
acts as an isometry of $(T\xtil,d_{T\xtil})$. For every $s,s'\in {_s S_i} (\xtil)$ we now set
$$
d_S (s,s')=\sup_{x\in T\Delta_i} d_{T\xtil} (ds (x),ds' (x)).
$$
It is well-known that $d_S$ induces on ${_s S_i} (\xtil)$ the $C^1$--topology.

Let now $s_0\in {_s S_i} (\xtil)$ and $\vare >0$ be fixed. By Lemma~\ref{contact:lemma}, the set
$\overline{F}\cdot s_0\subseteq {_s S_i} (\xtil)$ is compact. 
Since $\varphi$ is continuous, this 
easily implies that there exists $\eta>0$ such that $|\varphi (s_1)-\varphi (s_2)| \leq \vare$
for every $s_1\in \overline{F}\cdot s_0$, $s_2\in B_{d_S} (s_1,\eta)$, where 
$B_{d_S}(s_1,\eta)$ is the open ball
of radius $\eta$ centered at $s_1$. Take now $s\in B_{d_S} (s_0,\eta)$.
Since $G$ acts isometrically on ${_s S_i (X)}$, for every $g\in F$ we have $d_S (g\cdot s_0, g\cdot s)=
d_S (s_0,s)<\eta$, so $|\varphi (g\cdot s)-\varphi (g\cdot s_0)| \leq \vare$.
Together with the fact that $\mu_G (F)=1$, this readily implies 
$$
|\tr^i (\varphi) (s) -\tr^i (\varphi) (s_0)|  
 \leq 
 \int_F |\varphi (g\cdot s)-\varphi (g\cdot s_0) |\, d\mu_G (g)  \leq \vare .
 $$
We have thus proved that $\tr^i (\varphi)$ is continuous.


Now, if $\varphi$ is $G$--invariant then by the very definition we have $\tr^i (\varphi)=\varphi$.
The fact that $\tr^i (\varphi)$ is $G$--invariant if $\varphi$ is $\Gamma$--invariant
follows from the very same computations described
in~\cite[Subsection 6.3]{Bucher}.
\qed\smallskip 

Proposition~\ref{funziona:prop} provides a well-defined map 
$\tr^\ast\colon  {_s \ccst} (\xtil)^\Gamma\to {_s \ccst} (\xtil)^G$. It is readily seen
that $\tr^\ast$ is a chain map, and we still denote by
$\tr^\ast\colon {_s \hbst}(\xtil)^\Gamma\to {_s \hbst}(\xtil)^G$ the resulting map in cohomology.
Since $\tr^\ast$ restricts to the 
identity on $G$--invariant cochains, we have the following commutative diagram:

$$
\xymatrix{
{_s \hcst} (\xtil)^G \ar[r]_{\res^\ast} \ar@/^1.5pc/[rr]^{\rm Id} &  
{_s \hcst} (\xtil)^\Gamma 
\ar[r]_{\tr^\ast} \ar@{<->}[d] &  {_s \hcst} (\xtil)^G\\
& {_s \hcst} (X) &
}
$$
where the vertical row describes the isomorphism provided by Lemma~\ref{sollevosmooth:lemma}.
Since $\tr^\ast$ is obviously norm--decreasing, we get the following:

\begin{prop}\label{last:prop}
The map $\res^\ast \colon {_s \hcst}(\xtil)^G\to {_s \hcst}(\xtil)^\Gamma$ is an isometric embedding.
\qed
\end{prop}

\subsection{Proof of Theorem~\ref{newprop:teo}}
Since $G$ acts on $\xtil$ via orientation--preserving isometries, the cochain
${\rm Vol}_{\xtil}\in {_s C^n_c} (\xtil)$ is of course $G$--invariant, whence $\Gamma$--invariant.
Let us denote by $[{\rm Vol}_{\xtil}]_c^\Gamma$ (resp.~by $[{\rm Vol}_{\xtil}]_c^G$) the cohomology class
in ${_s \hc^n}(\xtil)^\Gamma$ (resp.~in ${_s \hc^n}(\xtil)^G$) represented by ${\rm Vol}_{\xtil}$.
Of course $\res^n ([{\rm Vol}_{\xtil}]_c^G)=[{\rm Vol}_{\xtil}]_c^\Gamma$, while since
$p\colon\xtil\to X$ is a local isometry we have  $p^\ast ({\rm Vol}_X)={\rm Vol}_{\xtil}$.
Therefore Proposition~\ref{last:prop} and
Lemma~\ref{sollevosmooth:lemma} imply
$$
||[{\rm Vol}_X]_c ||=||[{\rm Vol}_{\xtil}]_c^\Gamma||
=||\res^n ([{\rm Vol}_{\xtil}]_c^G)||= ||[{\rm Vol}_{\xtil}]_c^G||,
$$
whence Theorem~\ref{newprop:teo}, 
since $||[{\rm Vol}_{\xtil}]_c^G||$ only depends on the isometry type of $\xtil$.

\appendix
\section{Compact--open and $C^1$--topology}\label{basi:app}

\subsection{Compact--open topology}
Recall that if $X,Y$ are topological spaces, the \emph{compact--open} topology
on the space $F(X,Y)=\{f\colon X\to Y,\ f\ {\rm continuous}\}$ 
admits as a subbasis the set 
$$
\{\Omega (K,U),\ K\subset X\ {\rm compact},\ U\subset Y\ {\rm open}\},
$$
where 
$$
\Omega (K,U)=\{f\in F(X,Y)\, | \, f(K)\subseteq U\}.
$$

In this subsection, all the function spaces involved are endowed with the compact--open topology.
The following result is proved in~\cite[page 259]{Dug}:

\begin{lemma}\label{basi:top:lemma1}
Let $X,Y,Z$ be topological spaces, and $f\colon Y\to Z$, $g\colon X\to Y$ be continuous.
The maps $f_\ast \colon F(X,Y)\to F(X,Z)$, $g^\ast \colon F(Y,Z)\to F(X,Z)$ defined by $f_\ast (h)=f\circ h$,
$g^\ast (h)=h\circ g$ are continuous. 
\end{lemma}


\begin{lemma}\label{basi:top:lemma2}
Suppose $X$ is compact and Hausdorff, 
let $C\subseteq X$ be closed and set $F_C (X,Y)=\{h\in F(X,Y)\, |\, h|_C\ {\rm is\ constant}\}$.
Let $\pi\colon X\to X/C$ be the canonical projection, and
for $h\in F_C (X,Y)$, let $\psi (h)\in F(X/C,Y)$ be the unique map such that
$\psi (h)\circ \pi=h$. Then $\psi\colon F_C (X,Y)\to F(X/C,Y)$ is well-defined and continuous.
\end{lemma}
\noindent {\sc Proof:}
The fact that $\psi$ is well-defined is an immediate consequence of the definition
of quotient topology. Moreover, if $K\subseteq X/C$ is compact and $U\subseteq Y$ is open,
then $\psi^{-1} (\Omega (K,U))=\Omega (\pi^{-1}(K),U)$, which is open: in fact,
since $X$ is compact Hausdorff and $C$ is closed, $X/C$ is Hausdorff, so $K$ is closed, and
$\pi^{-1}(K)$ is 
closed in a compact space, whence compact.
\qed\smallskip

As in the proof of Proposition~\ref{strong1:prop},
let now 
$e_0^n,\ldots,e_n^n$ be the vertices of the standard simplex $\Delta_n$, 
let $Q_0^n$ be the face of $\Delta_n$ opposite to $e_0^n$, and let 
$r_n\colon Q^{n+1}_0\to \Delta_n$ be defined by
$r_n (t_1 e^{n+1}_1+\ldots t_{n+1} e^{n+1}_{n+1})=t_1 e^n_0+\ldots t_{n+1} e^n_n$.

\begin{lemma}\label{basi:top:lemma3}
Let $X,Y,Z$ be topological spaces, let $x_0\in X$ be fixed
and suppose 
$H\colon X\times [0,1]\to X$ is a continuous map such that $H(x,1)=x_0$ for every
$x\in X$. 
\begin{enumerate}
\item
The map $(h\times {\rm Id} \colon F(X,Y)\to F(X\times [0,1],Y\times [0,1])$
is continuous (where products are endowed with the product topology).
\item
Let $T_n\colon F(\Delta_n,X)\to F(\Delta_{n+1},X)$ be defined as follows: if $s\in F(\Delta_n,X)$ and
$p=t e^{n+1}_0+(1-t)q \in\Delta_{n+1}$, where $q\in Q^{n+1}_0$, then 
$$(T_n(s)) (p)=H(s(r_n(q)), t).$$
Then $T_n$ is well-defined and continuous.
\end{enumerate}
\end{lemma}
\noindent{\sc Proof:} 
(1): By~\cite[page 264]{Dug}, 
the compact--open topology of $F(X\times [0,1],Y\times [0,1])$ admits
as a subbasis the set 
$$
\mathcal{B}=\{\Omega (K,U\times U'),\ K\subseteq X\times [0,1]\ {\rm compact},\ U\subseteq Y\ {\rm open},\, 
U'\subseteq [0,1]\ {\rm open}\}.
$$ 
Now if $\pi_1\colon X\times [0,1]\to X$, $\pi_2\colon X\times [0,1]\to [0,1]$
are the natural projections and $\Omega (K,U \times U')\in\mathcal{B}$, then
$(h\times {\rm Id})^{-1} (\Omega (K,U\times U'))=\Omega (\pi_1 (K), U)$ if $\pi_2 (K)\subseteq U'$,
and $(h\times {\rm Id})^{-1} (\Omega (K,U\times U'))=\emptyset$ otherwise. In any case, $(h\times {\rm Id})^{-1} (\Omega (K,U\times U'))$
is open in $F(X,Y)$.

(2): Since $r_n^{-1}\colon \Delta_n\to Q_0^{n+1}\subset \Delta_{n+1}$ is a homeomorphism,
the map $\varphi\colon \Delta_n\times [0,1]\to \Delta_{n+1}$ defined by
$\varphi (q,t)=(1-t)r_n^{-1} (q)+t e_0^{n+1}$ is well-defined and continuous. Moreover, we have
$\varphi(q,1)=\varphi (q',1)$ for every
$q,q'\in \Delta_n$, so $\varphi$ induces a map 
$\overline{\varphi}\colon (\Delta_n\times [0,1])/(\Delta_n\times \{1\})\to\Delta_{n+1}$
which is easily shown to be bijective. Moreover, $\overline{\varphi}$ is continuous by the very definition
of quotient topology, and is closed since it is defined on a compact space with values in a Hausdorff
space. Thus, $\overline{\varphi}$ is a homeomorphism. Now, if $s\in F(\Delta_n,X)$, since $H(q,1)=x_0$
for every $q\in X$, the map
$H\circ (s\times {\rm Id})\colon \Delta_n\times [0,1]\to X$ defines a continuous
map $\overline{H\circ (s\times {\rm Id})}\colon (\Delta_n\times [0,1])/(\Delta_n\times\{1\})\to X$, and by
construction we have $T_n (s)= \overline{H\circ (s\times {\rm Id})}\circ \overline{\varphi}^{-1}$.
This shows that $T_n (s)$ is indeed well-defined and continuous.
 
Let us show now  
that $T_n (s)$ continuously depends on $s$. 
Since $\overline{\varphi}^{-1}$ is a homeomorphism, it is sufficient to show that
the map $\overline{H\circ (s\times {\rm Id})}$
continuously depends on $s$. But this is a consequence of Lemma~\ref{basi:top:lemma1},
Lemma~\ref{basi:top:lemma2} 
and point~(1).
\qed\smallskip

Let now $p\colon\xtil\to X$ be a covering, and denote by $p_\ast\colon C_n (\xtil)\to C_n (X)$
the induced map on singular chains. The following result was proved in~\cite{Loh}
under the hypothesis that $\xtil$ is metrizable. 

\begin{lemma}\label{open:lemma}
The restriction $p_\ast |_{S_n (\xtil)}\colon S_n (\xtil)\to S_n (X)$
is a covering map. 
\end{lemma}
\noindent {\sc Proof:} 
Let $s_0\in S_n (\xtil)$ be a simplex, and set $s_0(e_0)=x_0$. Since $p$ is a covering,
there exists an open neighbourhood $U_0$ of $x_0\in X$ 
such that $p^{-1} (U_0)=\bigsqcup_{j\in J} \widetilde{U}_0^j$, $\widetilde{U}_0^j$
is open in $\widetilde{X}$   
and $p|_{\widetilde{U}_0^j}\colon \widetilde{U}_0^j\to U_0$ is a homeomorphism for every $j\in J$.

We set $V_0=\{s\in S_n (X)\, | \, s(e_0)\in U_0\}$, $\widetilde{V}_0^j=\{\widetilde{s}\in S_n(\xtil)\, | \,
\widetilde{s} (e_0)\in \widetilde{U}_0^j\}$. Of course, $V_0$ and $\widetilde{V}_0^j$ are open
subsets of $S_n (X)$ and $S_n (\xtil)$ respectively. Moreover,
since the standard simplex is path connected and simply connected, for every $j\in J$ and 
every simplex in $s\in V_0$ there exists a unique lift  $\widetilde{s}^j\in \widetilde{V}_0^j$. 
This readily implies
that $p_\ast^{-1} (V_0)=\bigsqcup_{j\in J} \widetilde{V}_0^j$, and that 
$p_\ast |_{\widetilde{V}_0^j}\colon \widetilde{V}_0^j\to V_0$ is bijective for every
$j\in J$. Moreover, by Lemma~\ref{basi:top:lemma1} the map $p_\ast |_{S_n(\xtil)}$ is continuous,
so in order to conclude we are only left to show that 
$p_\ast |_{\widetilde{V}_0^j}\colon \widetilde{V}_0^j\to V_0$
is open for every $j\in J$.

We fix $\overline{j}\in J$, and denote $\widetilde{V}_0^{\overline{j}}$ 
(resp.~$\widetilde{U}_0^{\overline{j}}$)
simply by $\widetilde{V}_0$ (resp.~$\widetilde{U}_0$).
Since $p_\ast |_{\widetilde{V}_0}$ is injective, it preserves unions
and intersections, thus it is sufficent to prove that if $\widetilde{s}\in \widetilde{V}_0$ is a simplex and
$\widetilde{s}\in \Omega (K,Y)$, where $K\subseteq \Delta_n$
is compact and $Y\subseteq \xtil$ is open, then $p_\ast (\Omega(K,Y)\cap \widetilde{V}_0)$ is a neighbourhood
of $s=p_\ast (\widetilde{s})$.
Since $\widetilde{s}(\Delta_n)$ is compact, there exists a finite open cover $\{\widetilde{Z}_i\}_{i=0}^l$
of $\widetilde{s}(\Delta_n)$ such that 
$\widetilde{Z}_i$ homeomorphically projects onto an open
subset $Z_i\subseteq X$ for every $i$.
Moreover,
it is easily seen that there exists a decomposition $\Delta_n=\bigcup_{i=0}^N \Delta^i$ with the following
properties:
\begin{enumerate}
\item
each $\Delta^i$ is closed in $\Delta_n$;
\item
$e_0\in \Delta^i$ if and only if $i=0$;
\item
$\widetilde{s}(\Delta^0)\subseteq \widetilde{U}_0$;
\item
for every $i=0,\ldots,N$ there exists $j_i\in\{0,\ldots,l\}$ such that
$\widetilde{s}(\Delta^i)\subseteq \widetilde{Z}_{j_i}$.
\end{enumerate} 
Let $I=\{0,\ldots,N\}$. If $L\subseteq I$, we set $\Delta^L=\bigcap_{i\in L} \Delta^i$ and
$\widetilde{Z}_L=\bigcap_{i\in L} \widetilde{Z}_{j_i}$.
We also set 
$$
\widetilde{H}=\bigcap_{L\subseteq I} \Omega({\Delta^L},\widetilde{Z}_{L}),\quad 
\widetilde{R}=\bigcap_{L\subseteq I} \Omega ({K\cap \Delta^L},Y\cap\widetilde{Z}_{L}), \quad
\widetilde{W}_0=\Omega (\Delta^0,\widetilde{U}_0\cap\widetilde{Z}_{j_0}),
$$
and
$$
H= \bigcap_{L\subseteq I} \Omega({\Delta^L},p(\widetilde{Z}_L)),\quad 
R= \bigcap_{L\subseteq I} \Omega({K\cap\Delta^L},{p(Y\cap \widetilde{Z}_{L}})),\quad
W_0= \Omega (\Delta^0,p(\widetilde{U}_0\cap \widetilde{Z}_{j_0})). 
$$ 
Note that since $p$ is open, the sets $\widetilde{H},\widetilde{R},\widetilde{W}_0,H,R,W_0$ are open.
Moreover, by construction we have $\widetilde{s}\in \widetilde{W}_0\cap \widetilde{H}\cap\widetilde{R}\subseteq
\widetilde{V}_0\cap \Omega(K,Y)$ and $s\in W_0\cap H\cap R$.
Thus in order to conclude we only need to prove that
$$
p_\ast (\widetilde{W}_0\cap \widetilde{H}\cap \widetilde{R})= W_0\cap H\cap R.
$$
The inclusion $\subseteq$ is obvious. Let $s_1\in W_0\cap H\cap R$, and let
$\widetilde{s}_1\in \widetilde{V}_0$ be the unique lift of $s_1$ whose first
vertex lies in $\widetilde{U}_0$.

We now try to reconstruct $\widetilde{s}_1$ as explicitely as possible.
For every $i\in I$, 
we denote by $r_i\colon Z_{j_i}\to \widetilde{Z}_{j_i}$ the local inverse
of $p$. 
For every $q\in\Delta_n$, let $I(q)=\{i\in I\, |\, q\in\Delta^i\}$.
We claim that for every $q\in\Delta_n$ and $i,i'\in I(q)$ we have
$r_i (s_1(q))=r_{i'} (s_1 (q))$: 
in fact, since $s_1\in H$ 
we have $s_1 (q)\in p(\widetilde{Z}_{I(q)})$, so 
a point $a\in \widetilde{Z}_{j_i}\cap\widetilde{Z}_{j_{i'}}$ exists such that
$p(a)=s_1 (q)$. But $p(r_i (s_1(q)))=p(r_{i'}(s_1(q)))=s_1 (q)$, and $p|_{\widetilde{Z}_i}$,
$p|_{\widetilde{Z}_{i'}}$ are injective, so $r_i(s_1 (q))=a=r_{i'} (s_1 (q))$. 

Therefore, a well defined map $\overline{s}_1\colon \Delta_n\to \xtil$ exists such that
$\overline{s}_1 |_{\Delta^i}=r_i\circ s_1 |_{\Delta^i}$ for every $i\in I$. Since the
$\Delta^i$'s provide a finite closed cover of $\Delta_n$, this map is continuous. Moreover,
$p\circ\overline{s}_1=s_1$, and since $s_1\in W_0$  
we have $\overline{s}_1 (e_0)\in p^{-1}(p(\widetilde{U}_0\cap \widetilde{Z}_{j_0}))$.
Since by construction $\overline{s}_1 (e_0)\in\widetilde{Z}_{j_0}$ and $p|_{\widetilde{Z}_{j_0}}$ is 
injective, this implies
$\overline{s}_1 (e_0)\in\widetilde{U}_0$, so $\overline{s}_1 (e_0)=\widetilde{s}_1 (e_0)$.
This, together with the fact that $p\circ\overline{s}_1=p\circ\widetilde{s}_1$ and the
fact that $\Delta_n$ is path connected, implies that $\widetilde{s}_1=\overline{s}_1$.

We are then left to show that $\overline{s}_1$ indeed belongs to $\widetilde{W}_0\cap\widetilde{H}\cap
\widetilde{R}$. The fact that $\overline{s}_1$ belongs to $\widetilde{W}_0\cap\widetilde{H}$
is an immediate consequence of our construction. Now, if $q\in\Delta^L\cap K$, since
${s}_1\in R$ there exists $a\in Y\cap \widetilde{Z}_L$ such that $p(a)=s_1 (q)$.  
Moreover, by construction we have $p(a)=p(\overline{s}_1 (q))$ and $\overline{s}_1 (q)\in \widetilde{Z}_L$, so
$a=\overline{s}_1 (q)$ since $p|_{\widetilde{Z}_L}$ is injective. Therefore we have
$\overline{s}_1 (q)=a\in Y\cap \widetilde{Z}_L$, and this proves that
$\overline{s}_1$ belongs to $\widetilde{R}$, whence the conclusion.
\qed\smallskip

\subsection{$C^1$--topology on maps between smooth manifolds}
If $X,Y$ are smooth manifolds, we denote by $F_s (X,Y)$ the set
of smooth maps from $X$ to $Y$. 
If $TX,TY$ are the tangent bundles of $X,Y$ respectively,
the usual differential defines a map
$d\colon F_s(X,Y)\to F_s(TX,TY)$, which is obviously injective.
The \emph{$C^1$--topology} on $F_s (X,Y)$ is the pull-back via $d$ of the compact--open
topology of $F_s (TX,TY)$.
  
Putting together Lemma~\ref{basi:top:lemma1} and
the fact that the differential of the composition of smooth maps
is the composition of their differentials  
it is  not difficult to get the following:

\begin{lemma}\label{smoothbasi:lemma1}
Let $X,Y,Z$ be smooth manifolds, let $f\colon Y\to Z$, $g\colon X\to Y$ be smooth, and endow $F_s (X,Y)$, $F_s(Y,Z)$, 
$F_s(X,Z)$ with the $C^1$--topology. 
The maps $f_\ast \colon F(X,Y)\to F(X,Z)$, $g^\ast \colon F(Y,Z)\to F(X,Z)$ defined by $f_\ast (h)=f\circ h$,
$g^\ast (h)=h\circ g$ are continuous. 
\end{lemma}

In particular, if $h\colon X\to X$ is a smooth map and $\varphi\in {_s C^\ast_c (X)}$,
then $h^\ast (\varphi)\in {_s C^\ast_c} (X)$. 
Therefore, continuous smooth cohomology 
is a functor from the category of manifolds and smooth maps
to the category of graded $\R$--vector spaces and linear maps. Moreover, 
using Lemma~\ref{smoothbasi:lemma1}
and arguing just as in the proof of Proposition~\ref{hom:prop}
it is not difficult to prove the following:

\begin{lemma}\label{invhomsmooth}
Let $f,g\colon X\to Y$ be smooth maps between smooth manifolds, and 
let $f_c^\ast, g_c^\ast\colon {_s \hcst (Y)}\to {_s \hcst (X)}$
be the induced maps
in smooth cohomology.
If
$f$ is smoothly homotopic to $g$, then $f_c^\ast=g_c^\ast$.
\end{lemma} 

Moreover, if $G$ acts on $X$
as a group of diffeomorphisms, then it makes sense to define $G$--invariant continuous
smooth cochains on $X$.

Let now $X,\xtil$ be smooth manifolds, and
suppose $p\colon\xtil\to X$ is a smooth covering (\emph{i.e.}~a covering
which is also a local diffeomorphism).
We endow ${_s S_n} (X)$, ${_s S_n} (\xtil)$ with the $C^1$--topology.

\begin{lemma}\label{opensmooth:lemma}
The map $p_\ast |_{{_s S_n} (\xtil)}\colon {_s S_n} (\xtil)\to {_s S_n} (X)$
is a covering map. 
\end{lemma}
\noindent {\sc Proof:} 
If $X,Y$ are smooth manifolds, we say that
$g\colon TX\to TY$ is \emph{integrable} if there exists a smooth
$f\colon X\to Y$ such that $g=df$. Of course, a smooth map is integrable if and only
if it is locally integrable. 
Let us denote by $(dp)_\ast\colon F_s (T\Delta_n,T\xtil)\to F_s (T\Delta_n,T X)$
the composition with $dp$. 
We have the commutative diagram
$$
\xymatrix{
F_s (\Delta_n,\xtil) \ar[r]^d \ar[d]^{p_\ast} & F_s (T\Delta_n, T\xtil) \ar[d]^{(dp)_\ast}\\
F_s (\Delta_n,X) \ar[r]^d & F_s (T\Delta_n, TX)
}
$$
It is easily seen that, since $p$ is a smooth covering, $dp\colon T\xtil\to TX$
is also a smooth covering. Since $T\Delta_n$ is simply connected, 
a slight modification of 
the proof of Lemma~\ref{open:lemma}
shows that
$(dp)_\ast$ is a covering, where we are endowing $F_s (T\Delta_n,T\xtil)$, $F_s (T\Delta_n,T X)$
with the compact--open topology. Moreover, it is readily seen that $f\in F_s (T\Delta_n,T\xtil)$
is integrable if and only if $(dp)_\ast (f)$ is integrable. 
Since the subset of integrable maps in $F_s (T\Delta_n, T\xtil)$ (resp.~in $F_s (T\Delta_n, TX)$)
coincide by definition with $d(F_s (\Delta_n,\xtil))$ (resp.~with $d(F_s (\Delta_n,X))$),
$(dp)_\ast$ restricts to a covering $d(F_s (\Delta_n,\xtil))\to d(F_s (\Delta_n,X))$. The conclusion
follows since the horizontal rows of the diagram are by definition homeomorphisms on their
images.
\qed\smallskip

\bibliographystyle{amsalpha}
\bibliography{biblio}

\end{document}